\documentclass[11pt]{amsart}

\usepackage{amsmath,amsfonts,amssymb, amsthm,
latexsym,enumerate,graphics}
\usepackage{times}

\usepackage{xy}
\xyoption{all}
\CompileMatrices

\newcommand{\be}{\begin{equation}}
\newcommand{\ee}{\end{equation}}

\def\C{\ensuremath{\mathbb{C}}}
\def\D{\ensuremath{\mathbb{D}}}

\def\L{\ensuremath{\mathbb{L}}}

\def\P{\ensuremath{\mathbb{P}}}
\def\Q{\ensuremath{\mathbb{Q}}}
\def\R{\ensuremath{\mathbb{R}}}
\def\Z{\ensuremath{\mathbb{Z}}}

\def\k{\ensuremath{\underline{k}}}
\def\m{\ensuremath{\underline{m}}}

\def\Aut{\mathop{\mathrm{Aut}}\nolimits}
\def\age{\mathop{\mathrm{age}}\nolimits}
\def\Bl{\mathop{\mathrm Bl}\nolimits}

\def\rk{\mathop{\mathrm{rk}}}

\def\Spec{\mathop{\mathrm{Spec}}}

\def\Sym{\mathop{\mathrm{Sym}}\nolimits}

\def\({\left(}
\def\){\right)}

\newcommand\abs[1]{\lvert#1\rvert}
\newcommand\floor[1]{\lfloor#1\rfloor}
\newcommand\ceil[1]{\lceil#1\rceil}
\newcommand\fract[1]{\langle#1\rangle}
\newcommand\restr[2]{\left.#1\right|_{#2}}

\newcommand\stv[2]{\left.\kern-\nulldelimiterspace
                \left\{#1\vphantom{#2}\,\right|#2\right\}}

\newenvironment{Prf}{\textit{Proof.}\/}{\hfill$\Box$}

\newtheorem*{Thm-nonumber}{Theorem}
\newtheorem{Thm}{Theorem}[subsection]
\newtheorem{Def}[Thm]{Definition}

\newtheorem{Prop}[Thm]{Proposition}
\newtheorem{Lem}[Thm]{Lemma}
\newtheorem{Cor}[Thm]{Corollary}

\numberwithin{equation}{subsection}

\def\AA{\ensuremath{\mathcal A}}
\def\BB{\ensuremath{\mathcal B}}
\def\CC{\ensuremath{\mathcal C}}
\def\DD{\ensuremath{\mathcal D}}
\def\EE{\ensuremath{\mathcal E}}
\def\FF{\ensuremath{\mathcal F}}

\def\KK{\ensuremath{\mathcal K}}
\def\LL{\ensuremath{\mathcal L}}

\def\OO{\ensuremath{\mathcal O}}
\def\PP{\ensuremath{\mathcal P}}

\def\TT{\ensuremath{\mathcal T}}
\def\XX{\ensuremath{\mathcal X}}
\def\YY{\ensuremath{\mathcal Y}}

\def\ss{\ensuremath{\mathfrak s}}

\def\XDr{X_{D, r}}
\def\Mbar{\overline{M}}
\def\Cbar{\overline{C}}

\def\CbarAA{\overline{C}^{\mu_r}_{0, \AA}}
\def\CbarBB{\overline{C}^{\mu_r}_{0, \BB}}
\def\MbarAA{\overline{M}^{\mu_r}_{0, \AA, \EE}}
\def\MbarBB{\overline{M}^{\mu_r}_{0, \BB, \EE}}

\newcommand\otn[1]{[#1]}

\begin{document}

\title{Quantum cohomology of $[\C^N/\mu_r]$}
\thanks{The second author was supported by NSF grant No. 0502170.}

\author{Arend Bayer}
\address{University of Utah, Department of Mathematics, 155 South 1400 East,
Room 233, Salt Lake City, UT 84112}
\email{bayer@math.utah.edu}
\author{Charles Cadman}
\address{University of British Columbia, Department of Mathematics, 1984 Mathematics Rd,
Vancouver, BC, Canada V6P 4M8}
\email{cadman@math.ubc.ca}
\date{\today}

\begin{abstract}
We give a construction of the moduli space of stable
maps to the classifying stack $B\mu_r$ of a cyclic group by a sequence of
$r$-th root constructions on $\Mbar_{0, n}$.
We prove a closed formula for the total Chern class
of $\mu_r$-eigenspaces of the Hodge bundle, and thus of the obstruction
bundle of the genus zero Gromov-Witten theory of stacks of the
form $[\C^N/\mu_r]$.

We deduce linear recursions for all genus-zero
Gromov-Witten invariants.
\end{abstract}

\maketitle

\section{Introduction}

This paper combines two notions of stable maps---twisted (\cite{AV99}) and
weighted (\cite{Has03})---to produce a formula for the genus $0$ Gromov-Witten
invariants of $[\C^N/\mu_r]$.  More precisely, we derive a formula for the
equivariant Euler class of the obstruction bundle on
$\Mbar_{0,n}([\C^N/\mu_r])$ as a pull-back of classes on $\Mbar_{0, n}$.
Our definition of weighted twisted stable maps
is ad hoc, applying only to genus $0$ maps with target $[\C^N/\mu_r]$.
Nevertheless, there is a notion of
obstruction bundle on each of these spaces, and we have a wall-crossing
formula relating their equivariant Euler classes.  When all except one of the
weights is small, the equivariant Euler class is easy to compute, so we deduce
in this way our formula on $\Mbar_{0,n}([\C^N/\mu_r])$.

The starting point of this work is our theorem \ref{thm:modspace}, which provides an explicit construction of
$\Mbar_{0,n}([\C^N/\mu_r])$ from $\Mbar_{0,n}$ via root constructions.  This
motivates the generalization to weighted stable maps.  They are defined by
applying root constructions to the space of weighted stable curves.

\subsection{Introduction to $[\C^N/\mu_r]$}
The quantum cohomology of stack quotients of the form $[\C^N/G]$  has
recently seen a lot of interest due to the crepant resolution conjecture (see
\cite{Bryan-Graber:crepant-conj}). However, they are
also natural objects of study by themselves: whenever a smooth
$N$-dimensional stack $X$ has a local
orbifold chart $[U/G]$ where $G$ acts with an isolated fixed point, part of
the quantum cohomology of $X$ will be described by the quantum cohomology
of $[\C^N/G]$.  Moreover, if $X$ has a torus action which restricts to the natural torus action on the chart $[\C^N/G]$, then the equivariant quantum cohomology of $[\C^N/G]$ is relevant for computing the quantum cohomology of $X$ via localization.

The Chen-Ruan orbifold cohomology of $[\C^N/\mu_r]$ has a natural basis
$h_e$ for $e \in \mu_r$. Consider the moduli space
$\Mbar_{0, n}(e_1, \dots, e_n; B\mu_r)$ of twisted stable maps to
the origin $B\mu_r \cong [0/\mu_r] \subset [\C^N/\mu_r]$ in the sense
of \cite{AV99}, where the branching behavior
at the $i$-th section is prescribed by $e_i \in \mu_r$.
The non-trivial equivariant Gromov-Witten invariants of $[\C^N/\mu_r]$ are
given by integrals over these moduli spaces.
The normal bundle to the origin is $\C^N$, understood as a vector bundle
on $B\mu_r$ via the given $\mu_r$-action. If
we write the universal curve as
$\pi \colon \Cbar \to  \Mbar_{0, n}(e_1, \dots, e_n; B\mu_r)$,
and the universal map as $f \colon \Cbar \to B\mu_r$, the obstruction bundle
of the moduli space is $R^1 \pi_* f^* \C^N$. The Gromov-Witten invariant
for $h_{e_1}, \dots, h_{e_n}$ is typically\footnote{The formula
needs an additional factor in case there is a coordinate direction on
which every $e_i$ acts trivially.}
given by the integral of the equivariant Euler class
(with respect to the canonical action of the
$N$-dimensional torus $\TT$ on $\C^N$)
of the obstruction bundle
\begin{equation}			\label{eq:GWdef}
\langle h_{e_1} \otimes \dots \otimes h_{e_n} \rangle^{[\C^N/\mu_r]}
= \int_{\Mbar_{0, n}(e_1, \dots, e_n; B\mu_r)}
e_\TT ([R^1 \pi_* f^* \C^N]).
\end{equation}
We call these integrals \emph{generalized Hurwitz-Hodge integrals}, as the
obstruction bundle is a direct sum of $\mu_r$-eigenspaces of the dual of the
Hodge bundle, where the moduli space is to be understood as a compactification
of the Hurwitz space of $\mu_r$-covers of $\P^1$ by \emph{admissible covers}.

\subsection{Our methods and results}
The starting point of our work is the following explicit description
of this moduli space of stable maps to $B\mu_r$ via the $r$-th root
construction of \cite{chuck:using-stacks-to-impose}. Given a divisor
$D$ on a scheme $X$, the $r$-th root construction $X_{D, r}$ is a stack
over $X$ that is isomorphic to $X$ outside of $D$, but whose points over
$D$ are stacky with $\mu_r$ as automorphism group.
For every proper subset $T \subset [n-1]:=\{1, \dots, n-1\}$ having at least 2 elements, let $r_T$ be the order
$\prod_{i \in T} e_i$, and let $D^T\subset\Mbar_{0,n}$ be the divisor consisting of curves having a node which separates the marking labels $1,\ldots, n$ into $T$ and $[n]\setminus T$.

\begin{Thm-nonumber}
$\Mbar_{0, n}(e_1, \dots, e_n; B\mu_r)$ is a $\mu_r$-gerbe over the
stack constructed from $\Mbar_{0, n}$ by successively doing the $r_T$-th root
construction
at the boundary divisor $D^T \subset \Mbar_{0, n}$ for all proper subsets
$T \subset \otn{n-1}$ having at least 2 elements.
\end{Thm-nonumber}

We prove this in theorem \ref{thm:modspace}, and we also give
an explicit description of the universal curve and of the $r$-torsion line
bundle defining the morphism to $B\mu_r$; see definitions \ref{def:M(1)} and
\ref{def:M(2)}.  The root constructions along
the boundary divisor introduce the additional automorphisms of curves
with stacky nodes, called \emph{ghost automorphisms}.

Now assume that $\mu_r$ is acting linearly on $\C^N$ with weights
$w_1, \dots, w_N$.
To determine a formula for the Chern class of the obstruction bundle,
we use a reduction guided
by the notion of weighted stable curves in \cite{Has03} and
weighted stable maps in \cite{Mustatas:Chowring, Alexeev-Guy, weighted-GW}.
Weighted stable curves with $n$ marked points depend on weight data
$a_1, \dots, a_n$, and yield many birational models of the moduli space
$\Mbar_{0, n}$. Particular choices of weight data lead to
an explicit presentation of the moduli space
$\Mbar_{0, n}$ by a series of blow-ups starting with $\P^{n-3}$, such that
each intermediate blow-up step has an interpretation as a moduli space.

Motivated by this work and guided by theorem
\ref{thm:modspace}, we make an ad-hoc definition of a
``moduli space of weighted stable maps to $B\mu_r$''
in section \ref{subsect:weighted-stable-maps-to-Bmur}.
When the weights are chosen such that all fibers of the universal
curve are irreducible, the obstruction bundle can easily be computed from
general facts about the $r$-th root construction;
we do this in section \ref{subsect:Pn-3-computation} for the weight data
that gives a moduli space isomorphic to $\P^{n-3}$.

By a careful analysis of the wall--crossing for changing weights in section
\ref{sect:weight-change}, we can lift this to a closed formula for the
equivariant top Chern class in equation (\ref{eq:GWdef}) for the standard
(non-weighted) stable maps. We will now state this formula in the case of $N=1$
and the standard representation of $\mu_r$ with weight one: 

For $1\le i\le n$, let $\delta_i \in [0,1)$
be the age of $e_i$, i.e. 
$e^{2 \pi i \delta_i} = e_i$.
For all subsets $T \subset [n]$, let
$\delta_T = \sum_{i \in T} \delta_i$.  Let $\fract{x}$ denote the fractional
part of $x$ if $x$ is not an integer, and let $\fract{x}=1$ if $x$ is an
integer.  For $T \subset [n-1], \abs{T} \ge 2$, let $\psi_T$ be any class
in $H^*(\Mbar_{0,n})$ such that
the restriction $D^T \cdot \psi_T$ is the $\psi$-class of the node
over $D^T$ on the component corresponding to $T$; 
e.g. we can set
$\psi_T :=-\psi_n+\sum_{[n-1]\supsetneq S\supsetneq T} D^S.$
\begin{Thm-nonumber} 	
The equivariant Euler class of the obstruction
bundle for $[\C/\mu_r]$ is given as
\[
e_\TT \left([R^1 \pi_* f^* \C]\right)
= 
\prod_{p= \fract{\delta_{[n-1]}}}^{\delta_{[n-1]} - 1}
(t - p \psi_n) \cdot
\prod_{\substack{T \subsetneq \otn{n-1} \\ 2 \le \abs{T}}}
\prod_{p = \fract{\delta_T}}^{\delta_T - 1}
	\left(1 + \frac{ p D^T}{t + p \psi_T} \right)
\]
\end{Thm-nonumber}
The case of different weight follows by adjusting the ages in the
above formula, and the case of $\C^N$ by multiplying the individual 
classes in $H^*(\Mbar_{0,n})$; the full formula is given in
theorem \ref{mainthm}.

In the appendix, it is shown that this class can be expressed as a continuous, piece-wise analytic function from a real $(n-1)$-dimensional torus to
$H^*_T(\Mbar_{0, n}, \R)$ that encodes equivariant top Chern classes for all stacks
$[\C^N/\mu_r]$ where $N$, $n$ and the weights $w_1, \dots, w_N$ are fixed,
and $r$ and $e_1, \dots, e_n$ are arbitrary.  See the discussion after the proof of Lemma \ref{combinatorial_relation}.

By a generalized inclusion-exclusion principle, the Chern class
formula leads to linear recursions
for all Gromov-Witten invariants of $[\C^N/\mu_r]$ by a sum over partitions,
where every partition corresponds to a moduli space of comb curves.
They are particularly nice for local Calabi-Yau 3-folds $[\C^3/\mu_r]$.
We deduce an explicit formula for the non-equivariant invariants
of $[\C^3/\mu_3]$.  These invariants are the integrals in equation (\ref{eq:GWdef}) for which $N=r=3$, $n$ is a multiple of $3$, and all $e_i=e^{2\pi i/3}$.  The recursion we discovered for these numbers is
\begin{multline*}
\langle h_{\omega}^{\otimes n}
\rangle_{0,n}^{[\C^3/\mu_3]}
=
(-1)^{n+1}\((\frac{n-4}{3})!\)^3\frac{1}{3} \\
+ \sum_{p=1}^{\frac{n-3}3} \sum_{\m} \frac{(-1)^{\abs{\m}+1}}{|\Aut \m|}
	\prod_{j=1}^k \((m_j - \frac 23)!\)^3
	M(n-1, \m)
	\langle h_{\omega}^{\otimes n - 3p}
	\rangle_{0,n}^{[\C^3/\mu_3]},
\end{multline*}
where the second sum is over all partitions $\m = (m_1, \dots, m_k)$ of $p$, 
$M (n-1, \m)$ is the multinomial coefficient
\[
M(n-1, \m) =
	\binom{n-1}{3m_1 + 1, \dots, 3 m_k + 1, n-1-\sum_j (3m_j + 1)},
\]
and $x!=\prod_{p=\fract{x}}^x p$.  The base case is $\langle h_{\omega}^{\otimes 3}\rangle_{0,3}^{[\C^3/\mu_3]} = 1/3$.

\subsection{Relation to other work}
The construction of the moduli space and the universal curve via $r$-th roots
has been described locally by Abramovich in
\cite[\S 3.5]{abramovich:lectures-GW-orbifolds}.  The global description
along with the explicit description of the universal map to $B \mu_r$ seems
to be new.

Our approach is to describe maps to $B \mu_r$ by $r$-torsion line
bundles, which is particularly convenient in combination with the $r$-th
root constructions. Of course, one can instead work more geometrically
with cyclic covers, which is the point of view adopted in \cite{ACV,
Chuck-Renzo}.

The Gromov-Witten theory of $[\C^N/G]$ has recently generated a lot of
interest due to the ``crepant resolution conjecture''; we refer
to \cite{Bryan-Graber:crepant-conj} for an introduction to the conjecture
and overview of the existing literature. When $X$ is an orbifold that
admits a crepant resolution $\pi \colon Y \to X$, Ruan first conjectured
that the quantum cohomology rings of $Y$ and $X$ are isomorphic
\cite{Ruan:stringy_geometry}, and suggested that the $q$-variables for
$\pi$-exceptional divisor classes on $Y$ need to be specialized to -1
\cite{Ruan:crepant} to recover the orbifold cohomology ring of $X$. In
\cite{Bryan-Graber:crepant-conj}, the authors extended this conjecture:
their claim can be formulated as a local linear isomorphism between the
Frobenius manifolds of the quantum cohomology of $X$ and $Y$ (after analytic
continuation). This isomorphism does not respect the natural origins of the
two Frobenius manifolds, which corresponds to Ruan's specialization of
$q$-variables.  When the action of $G$ on $\C^N$ leaves the volume form
invariant, the stacks of the form $[\C^N/G]$ yield many non-trivial test cases
for the conjecture.  In this form it is only expected to hold for orbifolds
satisfying the strong Lefschetz theorem; a more general formulation can be found in \cite[section 5]{CCIT:wall-crossing}.

The results so far have been obtained by the use of one the following two
techniques: either a combination of localization computations and use of the
WDVV-equations, or by using Tseng's computation of the Chern \emph{character}
of the obstruction bundle and Givental's framework for
Gromov-Witten theory.

In \cite{BGP}, the authors explicitly determine the genus-zero
Gromov-Witten potential of $[\C^2/\mu_3]$ and verify the crepant resolution
conjecture. In \cite{Bryan-Graber:crepant-conj}, the case $[\C^2/\mu_2]$
is derived from the Hodge integral computations of
\cite{FaberPandharipande}.

More generally, the case of $A_n$-singularities $[\C^2/\mu_{n+1}]$
is shown in
\cite{CCIT:computing, CCIT:An} based on the Chern character computation.
Various other results have been announced in
\cite{Bryan-Graber:crepant-conj}.  For $[\C^3/\mu_3]$, part of the potential is
computed in \cite{CCIT:computing}, up to the problem of inversion of the
``mirror map''.  While their technique is completely
different to ours, our results are surprisingly close, as explained in section
\ref{twisted-I} and \ref{recursion-and-mirror-map}; our recursion can
be interpreted as a combinatorial inversion of the mirror map.

Different recursions for invariants of $[\C^3/\mu_3]$ have been established
by the second
author and R. Cavalieri in \cite{Chuck-Renzo}, using localization on the
space of twisted stable maps to $\mu_3$-gerbes over $\P^1$.

While our results are quite general, we make no attempt
at verifying the crepant resolution conjecture.

\subsection{Acknowledgements}
This project owes special thanks to R. Cavalieri; his explanations of
\cite{BGP} got us started on the direction of this work.
The first author would like to thank A.
Bertram, Y. Iwao, Y.-P. Lee and G. Todorov for
discussions about Givental's formalism and \cite{CCIT:computing} that
helped understanding the relation to our results.

\subsection{Notations and conventions}		\label{sect:notations}
We write $[n]$ for the set $[n] = \{1, 2, \dots, n\}$. We write
$\fract{x} = x - \ceil{x} + 1 \in (0,1]$
for the fractional part of $x$, set to 1 if $x$ is
integral. At various places
we will write $\prod_{p = \fract{x}}^x f(p)$ for the product
$\prod_{0 < p \le x, \fract{p} = \fract{x}} f(p)$.
If $x < 0$, the notation $\prod_{p = \fract{x}}^x f(p)$ means
$\prod_{p = x + 1}^{\fract{x}-1} \frac 1{f(p)}$ (which is consistent with
$\prod_{p = \fract{x}}^x f(p) = f(x) \cdot \prod_{p = \fract{x}}^{x-1} f(p)$
for all $x \in \R$).

For $x > 0$, we write $x!$ for the fractional factorial
$x! = \prod_{p = \fract{x}}^x p$.

We identify the rational Chow groups of the moduli stacks of twisted
stable maps $\Mbar_{0, n}(e_1, \dots, e_n; B\mu_r)$ with that of
its coarse moduli $\Mbar_{0, n}$ via pull-back, and similarly for all
other moduli stacks we construct.
In the appendix, we introduce and explain 
a few non-standard notations for divisors on $\Mbar_{0, n}$ that
are particularly well-suited for our setting; most of it is only used
in section \ref{subsect:formula-remarks}, the exception being
$$ \psi_T :=-\psi_n+\sum_{[n-1]\supsetneq S\supsetneq T} D^S$$
for any $T \subset [n-1]$.

\section{Moduli space of stable maps to $B\mu_r$ via $r$-th roots}
\label{mod-space-via-rth-roots}

In this section, we show how to construct a component of the moduli space of
genus zero stable maps to $B\mu_r$ from the moduli space $\Mbar_{0, n}$
of stable curves of genus zero by a series of $r$-th root constructions.

\subsection{User's guide to the $r$-th root construction}

Given an effective Cartier divisor $D$ of a Deligne-Mumford stack $X$, and a
positive integer $r$ which is invertible on $X$, the $r$-th root construction
of \cite{chuck:using-stacks-to-impose} produces a DM-stack
$\XDr$ with the following properties:

\begin{enumerate}
\item There is a canonical map $\pi \colon X_{D, r} \to X$ that is an
isomorphism over $X \setminus D$.
\item Every point in $\XDr$ lying over $D \subset X$ has stabilizer
$\mu_r$.
\item
The preimage of $D$ is an infinitesimal neighborhood of the
$\mu_r$-gerbe\footnote{A gerbe over $D$ is a stack $\DD$ over $D$ which \'etale
locally admits a section and has the property that any two local sections are
locally 2-isomorphic.  A gerbe $\DD\to D$ is a $\mu_r$-gerbe if $\mu_r$ acts as the
2-automorphism group of every section in a compatible way.} $\DD$ over
$D$ parameterizing $r$-th roots of the fibers of $\OO_X(D)|_D$: this is the
stack whose objects are triples $(f \colon S \to D, L, \phi)$ where
$f$ is a morphism, $L$ is a line bundle on $S$, and
$\phi \colon L^r \to f^* \OO_X(D)|_D$ is an isomorphism.
(Only when $\restr{\OO_X(D)}{D}$ is the $r$-th power of a line bundle is
$\DD$ isomorphic to $D \times B\mu_r$.)
\item \label{definingproperty}
On $\XDr$, there is a line bundle $\OO_{\XDr}(\DD)$ with a section $s_{\DD}$
and an isomorphism $\phi \colon \OO(\DD)^r \to \pi^*(\OO(D))$ such that
$\phi(s_{\DD}^r) = \pi^*(s_D)$. (Here $s_D \in \OO(D)$ is the tautological section
vanishing along $D$.)
\end{enumerate}

The universality of the data in (\ref{definingproperty}) is the defining
property: giving a morphism $f \colon S \to \XDr$ is equivalent to giving a
quadruple $(g, L, s, \phi)$ where $g = \pi \circ f$ is a morphism to $X$, $L$
is a line bundle on $S$, $s$ a section, and $\phi \colon L^r \to g^* \OO(D)$
is an isomorphism sending $s^r$ to $g^*(s_D)$.

Locally, when $X = \Spec A$ is affine and the divisor $D$ is given by an
equation $(x=0)$, the $r$-th root construction is given by
stack quotient $[\bigl(\Spec A[u]/(u^r-x)\bigr)/\mu_r]$ of the cyclic
$\mu_r$-cover branched at $D$, but of course globally such a cover may not
exist.

To the best of our knowledge, the $r$-th root construction is originally
due to A. Vistoli and spread as a rumor for quite some time. His notation is
$\sqrt[r]{(X, D)}$.

We call $\OO_{\XDr}(\DD)$ the tautological line bundle of the $r$-th root
construction at $D$. The zero stack $Z \subset \XDr$ of $s_{\DD}^r$ is the
preimage of $D \subset X$. The zero stack of $s_{\DD}$
is the gerbe $\DD$.
To simplify notation we write $\OO(\frac{1}{r}D)$ to refer to $\OO(\DD)$.  More
generally, if $d\in\frac{1}{r}\Z$, we write $\OO(dD)$ for $\OO(\DD)^{\otimes dr}$.
This notation is particularly nice to describe the push-forward of line
bundles along $\pi$:
\begin{equation}			\label{eq:push-forward}
\pi_* \OO(d D) = \OO(\floor{d} D)
\end{equation}
(This follows from \cite[Theorem 3.1.1]{chuck:using-stacks-to-impose}.)

If $X$ is an algebraic space, then the coarse  moduli space of $\XDr$ is $X$.
When $X$ is smooth and $D \subset X$ is
smooth, then $\XDr$ is smooth.  The construction commutes with base change
for a morphism $f \colon Y \to X$ such that $f^{-1}(D)$ is a Cartier divisor.
(The construction can be generalized a little to make it compatible with
arbitrary base change: see
$X_{(L, s, r)}$ in \cite{chuck:using-stacks-to-impose}.)

If $\D = (D_1, \dots, D_n)$ is an $n$-tuple of Cartier divisors and
$\vec{r}=(r_1,\ldots,r_n)$, we can iterate the root constructions to obtain
a stack denoted $X_{\D,\vec{r}}$.  This stack can also be realized as the
$n$-fold fiber product over $X$ of the root stacks $X_{D_i,r_i}$. If $X$
is smooth, each individual $D_i$ is smooth and the $D_i$ have normal
crossing, then $X_{\D, \vec{r}}$ is smooth, too.
If $d_i\in\frac{1}{r_i}\Z$, we extend the above notation by writing
$\OO(\sum_i d_iD_i)$ for the tensor product of the line bundles
$\OO(\DD_i)^{d_ir_i}$, where $\OO(\DD_i)$ is the tautological bundle
corresponding to the root construction along $D_i$.

It is instructive (and important for the construction of
$\Mbar_{0,n}(B\mu_r)$
later on) to compare the stacks $X_{\D,\vec{r}}$ and $\XDr$ in the case where
$n=2$, $r=r_1=r_2$, and $D=D_1\cup D_2$.  On $X_{\D, r}$, the line bundle
$\OO(\frac{1}{r}(D_1+D_2))$
with section $s_{\DD_1}\cdot s_{\DD_2}$ defines an $r$-th root
of $D$, and thus there is a natural map
\begin{equation}
\label{two_divisors}
X_{\D, \vec{r}} \to X_{D, r}.
\end{equation}
However, this is not an isomorphism if the divisors intersect, and one way
to see this is by looking at stabilizer groups of points of the two stacks.
If $x\in D_1\cap D_2$, then the stabilizer group of the point in
$X_{\D,\vec{r}}$ lying over $x$ is $\mu_r\times\mu_r$.  On the other hand,
the stabilizer group of any point in the preimage of $D$ in $\XDr$ is
$\mu_r$.  If $X$ is smooth and $D_1, D_2$ are smooth with normal crossings,
then these stacks are also distinguished by the fact that $X_{\D,\vec{r}}$
is smooth, while $\XDr$ is singular over $D_1\cap D_2$.

\subsection{Stable maps to $B\mu_r$}

Consider the moduli space $\Mbar_{0, n}(B\mu_r)$ of balanced twisted
stable maps of genus $0$
to $B\mu_r$ in the sense of \cite{AV99}, where we work over $\C$.
Such a map over a scheme $S$ can be described by the following data:
\begin{itemize}
\item A stacky nodal curve $\CC$ over $S$, with $n$ divisors
$\Sigma_1, \dots, \Sigma_n$ in the smooth locus of $\CC$, and
\item a line bundle $L$ on $\CC$ together with an isomorphism
$\phi \colon L^r \to \OO_\CC$.
\end{itemize}

These have to satisfy various properties:
\begin{itemize}
\item
The stacky curve
$\CC$ is a scheme away from its nodes and the divisors $\Sigma_i$,
and its nodes are \emph{balanced}.
\item
Each $\Sigma_i$ is a cyclotomic gerbe over $S$.
\item
If $C$ is the coarse moduli space of $\CC$, then the image of
every divisor $\Sigma_i \subset \CC$ in $C$
is isomorphic to the image
of a section $x_i \colon S \to C$, so that $(C, x_1, \dots, x_n)$ becomes
a stable curve of genus zero with $n$ marked points.
\item The map $\CC \to B\mu_r$ induced by $(L, \phi)$ is representable.
\end{itemize}

The line bundle $L$ is the pull-back of the line bundle on $B\mu_r$
given by the canonical one-dimensional representation of $\mu_r$.  Every
point $x\in\Sigma_i$ has an automorphism group isomorphic to $\mu_p$ for
some $p$ dividing $r$.  This identification is canonical if the
representation of $\mu_p$ corresponding to the fiber of the normal bundle
$\OO_{\Sigma_i}(\Sigma_i)$ at $x$ equals the standard representation.  Let
$\omega$ be the primitive $r$-th root of unity $\omega = e^{\frac{2\pi i}r}$.
Then $\omega^{r/p}$ acts on
the fiber $L_{x_i}$ as multiplication by $e_i$ for some
$e_i\in \mu_r$.  Equivalently, the map of stabilizer groups
$\mu_p\to\mu_r$ (which is injective
by representability of $\CC\to B\mu_r$) sends the canonical generator
of $\mu_p$ to $e_i$.  These group elements $e_1, \dots,
e_n$ are constant on every connected component of $\Mbar_{0, n}(B\mu_r)$.

 From now on, we assume we are given $e_1, \dots, e_n \in \mu_r$ and
restrict our attention to the connected component
$\Mbar_{0, n}(e_1, \dots, e_n; B\mu_r)$.  There is a natural map
$\Mbar_{0, n}(e_1, \dots, e_n; B\mu_r) \to \Mbar_{0, n}$ induced
by the coarse moduli space of the universal curve. Our theorem
will describe this map explicitly via a series of root
constructions.  One explanation for these root constructions is
that twisted curves
have ghost automorphisms for each twisted node (cf.
\cite[section 3.5]{abramovich:lectures-GW-orbifolds}).

\subsection{Construction of the moduli space via root constructions}

Let $r$ and $e_i \in \mu_r$, $i = 1,\ldots,r$ be given.
For convenience we allow $e_i = 1$, i.e. untwisted points; then the universal
curve is given by the forgetful morphism
\[
\pi_{0, (e_1, \dots, e_n)} \colon \Mbar_{0, n}(e_1, \dots, e_n, 1; B\mu_r)
\to \Mbar_{0, n}(e_1, \dots, e_n; B\mu_r).
\]
The component is empty unless $\prod_i e_i = 1$.

Consider the universal stable curve of genus zero
$\pi_{0,n} \colon \overline{M}_{0, n+1} \to \overline{M}_{0, n}$.
The boundary divisors of $\overline{M}_{0, n}$ are indexed by subsets
$T \subset \otn{n}$ such that $2 \le \abs{T} \le n-2$ and $n\not\in T$ (see
Appendix A).
For every such $T$, let $r_T$ be the order of
$\prod_{i \in T} e_i$.

\begin{Def} 				\label{def:M(1)}
Let $\Mbar^{(1)}$ be the stack constructed from $\Mbar_{0,n}$ by doing the
$r_T$-th root construction at every boundary divisor $D^T$.
We construct $\Cbar^{(1)}$ from $\Mbar_{0, n+1}$ in the same way after setting
$e_{n+1} = 1$.
\end{Def}
(In particular, we take the $r_i$-th root construction
at every section $s_i = D^{i, n+1}$ where $r_i$ is the
order of $e_i$.)

\begin{Lem}					\label{lem:pi1}
There is a canonical map $\pi^{(1)} \colon \Cbar^{(1)} \to \Mbar^{(1)}$.
\end{Lem}
\begin{Prf}
Equivalently,
we construct a map to the fiber product
\[
\Cbar^{(0)} = \Mbar_{0, n+1} \times_{\Mbar_{0,n}} \Mbar^{(1)};
\]
then $\Cbar^{(0)}$
will be the relative coarse moduli space.
Since $\pi_{0,n}^{-1} (D^T) = D^T \cup D^{T \cup \{n+1\}}$ and the
$r$-th root construction is compatible with such a base change, this
fiber product can be constructed from $\Mbar_{0, n+1}$ by the $r_T$-th root
constructions at all divisors $D^T \cup D^{T \cup \{n+1\}}$ for
$T \subset \otn{n}$. To construct $\Cbar^{(1)}$, we instead took the
$r_T$-th root construction at $D^T$ and $D^{T \cup \{n+1\}}$
separately, hence the morphism $\Cbar^{(1)} \to \Cbar^{(0)}$ is given by forgetting
the root construction along all sections, followed by a
composition of morphisms as in (\ref{two_divisors}) above.
\end{Prf}

Note that $\Cbar^{(1)}$ has additional automorphisms along the nodes of the
curves lying over $D^T$.  When we restrict this family to
the $\mu_{r_T}$-gerbe in $\Mbar^{(1)}$ lying over $D^T$, the fibers become stacky
curves with a twisted node.  The node is balanced because after base change
to a scheme over the base, the remaining automorphism group is the kernel
of the multiplication $\mu_{r_T} \times \mu_{r_T} \to \mu_{r_T}$, which acts
with opposite weight on the two branches.  The so--called ``ghost
automorphisms'' are accounted for by the additional automorphism introduced
in the moduli space. We have thus proved:

\begin{Prop}
The morphism $\pi^{(1)} \colon \Cbar^{(1)} \to \Mbar^{(1)}$ is a family of balanced
twisted curves.
\end{Prop}
Note that we adapted \cite[Definition 4.1.2]{AV99} to a family over
a Deligne-Mumford stack: all conditions have to be checked after \'etale
base change to a scheme covering $\Mbar^{(1)}$.

Each fiber of $\pi^{(1)}$ admits a morphism to $B\mu_r$ having the correct
restrictions to $\Sigma_i$ (given by $e_1,\ldots, e_n$).  However, these
morphisms do not in general glue to a morphism $\Cbar^{(1)}\to B\mu_r$.  They
will glue precisely when the $\mu_r$-gerbe of Definition \ref{def:M(2)} is trivial.

\subsection{The universal line bundle}

For $T \subset \otn{n}$, we will write always $T^C = \otn{n} \setminus T$
for its complement.

\begin{Lem} 				
Let $C$ be a geometric fiber of $\pi^{(1)}$ and let
$$\LL=\OO_{\Cbar^{(1)}}(\frac{1}{r_T}D^{T,n+1}).$$
\begin{enumerate}
\item If there is no node $x \in C$ such that one of the two connected
components of $C \setminus \{x\}$ contains exactly the markings
of $T$ (and the other those of $T^C$), then
$\restr{\LL}C$ is trivial.
\item
Otherwise,
let $C_3, C_4$ be the two connected components
of $C$ after normalization at $x$, such that
$C_3$ contains all the markings of $T$, and $C_4$ those of
$T^C$,
and let $C_1 \subset C_3$ and $C_2 \subset C_4$ be the two irreducible
components of $C$ meeting at $x$.
Then
\begin{eqnarray}
\restr{\LL}{C_1} & \cong & \OO_{C_1}(-\tfrac 1{r_T} x)
			\label{restrneg} \\
\restr{\LL}{C_2} & \cong & \OO_{C_2}(\tfrac 1{r_T} x)
						\label{restrpos} \\
\restr{\LL}{C'} & \cong & \OO_{C'}		  \label{restrivial}
\end{eqnarray}
\end{enumerate}
where $C'$ is any irreducible component of $C$ other than $C_1, C_2$.
\end{Lem}

\begin{Prf}
The first statement is obvious as $C$ does not meet the divisor
$D^{T}$ in that case.

In the second case, $C_3 = C \cap D^{T,n+1}$, and so
equation (\ref{restrpos}) is obvious, as is
(\ref{restrivial}) for all $C' \subset C_4$. The claim then follows
by symmetry and the fact that the restriction of
$\OO(\frac{1}{r_T}(D^T+D^{T,n+1}))$
to $C$ is trivial, as it is the pull-back
of the tautological line bundle $\OO(\frac{1}{r_T}D^T)$ on $\Mbar^{(1)}$.
\end{Prf}

Now choose $d_i \in \frac 1r \cdot \Z$ such that $e^{2 \pi i d_i} = e_i$
and $\sum_{i=1}^n d_i = 0$ (which is possible since
$\prod_{i=1}^n e_i = 1$).
For $T \subset \otn{n}$, let $d_T = \sum_{i \in T} d_i$.

\begin{Lem}
Define the line bundle $L_1$ on $\Cbar^{(1)}$ as
\begin{equation} 		\label{def:L1}
L_1 := \OO \Bigl(\sum_{i=1}^n d_is_i + \sum_{\substack{T \subset \otn{n}\\ 2 \le \abs{T}
\le n-2\\ n \not\in T}} d_T D^{T, n+1}\Bigr)
\end{equation}
Then ${L_1}^r$ is the pull-back of a line bundle on $\Mbar^{(1)}$:
\[
{L_1}^r = (\pi^{(1)})^*(L_2)
\]
\end{Lem}
We write $n \not\in T$ to stress that this is for now just an arbitrary
way to pick exactly one of $T, T^C$ for all subsets $T$.

\begin{Prf}
First note that $L_1^r$ is pulled back from the coarse moduli space
$\Mbar_{0,n+1}$
since $d_T\in\frac{1}{r}\Z$ for all $T$. A line bundle on a family of
nodal curves of genus zero over a smooth scheme is pulled back from the base
if and only if its degree on any irreducible component of every fiber is
zero (in which case it is the pull-back of its own push-forward to the base). Hence it is sufficient to
check that the degree of ${L_1}^r$ (or, equivalently, the degree of $L_1$)
is zero on any irreducible component
$C_0$ of any fiber $C$ of $\pi^1$
(in which case $L_1^r$ is even pulled back from $\Mbar_{0,n}$).

Let $x_1, \dots, x_m$ be the nodes of $C$ contained in $C_0$. Let
$T_i \subset \{1, \dots, n\}, 1 \le i \le m$ be the markings contained in
the irreducible components
which are connected to $C_0$ via the node $x_i$, and let
$T_0$ be the markings
contained in $C_0$. For every $j$ with $1 \le j \le m$ exactly one of
$\OO(d_{T_j}D^{T_j, n+1})$ and $\OO(d_{T_j^C} D^{T_j^C, n+1})$ will appear
in the right-hand side of equation (\ref{def:L1}) defining $L_1$;
by the previous lemma
and $d_{T} = - d_{T^C}$, both restrict to
$\OO(d_{T_j} x_j)$ on $C_0$.
By the same lemma, all other $\OO(d_TD^{T, n+1})$ restrict trivially to $C_0$.
Hence
\begin{flalign*}
\restr{L_1}{C_0}
& = \OO(\sum_{i\in T_0} d_is_i + \sum_{1\le j\le m} d_{T_j} x_j).
\end{flalign*}
This line bundle has degree
$\sum_{j=0}^m \sum_{i \in T_j} d_i = 0$ since $[n]$ is the disjoint union
of all $T_j$.
\end{Prf}

\subsection{Base change to the gerbe}

\begin{Def}			\label{def:M(2)}
Let $\Mbar^{(2)}$ be the $\mu_r$-gerbe over $\Mbar^{(1)}$ of
$r$-th roots of $L_2$.
Let $\pi^{(2)} \colon \Cbar^{(2)} \to \Mbar^{(2)}$
be the base change of $\pi^{(1)} \colon \Cbar^{(1)} \to \Mbar^{(1)}$
via $\Mbar^{(2)} \to \Mbar^{(1)}$, and
let $L_2^{1/r}$ be the universal line bundle that is an
$r$-th root of $L_2$.
\end{Def}
By abuse of notation, we write $L_1$ also for the pull-back
of $L_1$ to $\Cbar^{(2)}$.  The line bundle $L=L_1 \otimes {\pi^{(2)}}^* L_2^{-1/r}$ together with the obvious isomorphism $L^r\to\OO_{\Cbar^{(2)}}$ defines a morphism $\Cbar^{(2)} \to B\mu_r$.

\begin{Thm}				\label{thm:modspace}
The following diagram
$$\xymatrix{
\Cbar^{(2)} \ar[r] \ar[d]_{\pi^{(2)}} & B\mu_r \\
\Mbar^{(2)}}$$
is a family of twisted stable maps over $\Mbar^{(2)}$ which
defines an isomorphism $$m \colon \Mbar^{(2)} \to \Mbar_{0, n}(e_1, \dots, e_n; B\mu_r).$$
\end{Thm}

\begin{Prf}
We already showed that $\pi^{(1)}$ (and thus
$\pi^{(2)}$) is a family of balanced twisted curves.

The morphism $\Cbar^{(2)} \to \Mbar^{(2)} \times B\mu_r$
is representable: away from the sections and nodes, the map $\pi^{(2)}$
is already representable, and since all nodes and sections do not intersect each
other, we can treat them separately. At a section $s_i$, the relative
inertia group of $\pi^{(2)}$ is isomorphic to $B\mu_{r_i}$; since that
group acts faithfully on $L_1$ and thus on $L$, the map on inertia groups
is injective. A similar argument holds for all nodes.

We thus get a morphism $m$ as claimed in the theorem.
By lemma \ref{isom_criterion}, $m$ is an isomorphism if both
stacks are smooth and the morphism is birational and a bijection of
$\C$-valued points which induces isomorphisms of their stabilizer groups.
Since $\Mbar^{(2)}$ is \'etale over $\Mbar^{(1)}$, which is a root
construction on $\Mbar_{0,n}$ at smooth divisors with transversal
intersection, it follows
that $\Mbar^{(2)}$ is smooth.  The first order deformations of an $n$-marked,
genus $0$ twisted stable map to $B\mu_r$ are the same as those of the marked
twisted curve, which has dimension equal to $n-3$ (see
\cite[section 3]{ACV}).  As this equals the dimension
of $\Mbar_{0, n}(e_1, \dots, e_n; B\mu_r)$, it is also smooth.

For bijectivity, it suffices to show for each $n$-marked genus $0$ curve $C$,
there is a unique twisted stable map to $B\mu_r$ with coarse moduli space $C$
and contact types $e_1,\ldots,e_n$.  For $C$ irreducible, this uniqueness is
shown in \cite[2.1.5]{cadman-chen}.  Otherwise, one can show by induction that
the contact types at the nodes are uniquely determined by those at marked
points.  So the morphism is unique over each component of $C$, and it suffices
to show that it glues uniquely over the nodes.  Since the morphism to $B\mu_r$
is equivalent to a line bundle $L$ and a non-vanishing section of $L^r$, the
gluing is clearly unique up to isomorphism.

It remains only to check the map on automorphism groups induced by $m$.
If $x$ is a closed point of $\Mbar^{(2)}$,
the automorphism group $G_{m(x)}$ of $m(x)$ in the moduli stack
$\Mbar_{0, n}(e_1, \dots, e_n; B\mu_r)$ can be identified with the group of
$\mu_r$-automorphisms of the $\mu_r$-cover $\widetilde C_x$ of the fiber
$C_x$ of $\Cbar^{(2)}$ over $x$. If $S$ is the set of irreducible components
of $C_x$, this identifies $G_{m(x)}$ with the subgroup of
$\mu_r^{S}$ that acts compatibly over every node; since the preimage of
a node $n_T \in C_x$ in $\widetilde C_x$ is isomorphic to
$\mu_r/\mu_{r_T}$, we can identify $G_{m(x)}$
with the kernel of the map
\begin{equation*}	
\Sigma = \prod_T \Sigma_T
\colon \mu_r^S \to \prod_{T \subset [n-1] | x \in D^T} \mu_r/\mu_{r_T}
\end{equation*}
where $\Sigma_T$ is given by the quotient of the group elements corresponding
to the two irreducible components meeting in $n_T$.

The $\mu_r$-cover $\widetilde C_x$ is given by the $r$-th roots of unity
inside the line bundle $\restr{L}{C_x}$; hence to understand the map
$G_x \to G_{m(x)}$ it is
sufficient to look at how the automorphism group $G_x$ of $x$ acts on $L$.  By
the construction of $\Mbar^{(2)}$, the automorphism group of $x$ is
\[ G_x = \mu_r \times \prod_{T \subset [n-1]| x \in D^T} \mu_{r_T}.\]
By the definition of $L$, the first factor acts via $L_2$ and thus
diagonally, whereas $\mu_{r_T}$ acts diagonally on the
irreducible components of $D^{T, n+1} \cap C_x$ (and trivially on all others),
by its induced action on $\OO(\DD^{T, n+1})$.  Let $C_0$ be the component of $C_x$ which contains the $n$-th marking.

Before showing that $m:G_x\to G_{m(x)}$ is an isomorphism, we introduce some notation.  For $g\in G_x$, write $g_0$ for the projection of $g$ onto $\mu_r$, and write $g_T$ for its projection onto $\mu_{r_T}$.  For each irreducible component $C_i$ of $C_x$, let $C_{j_0},C_{j_1},C_{j_2},\ldots,C_{j_{k_i}}$ be the unique shortest path from $C_i$ to $C_0$.  That is to say, $j_0=i$, $j_{k_i}=0$, $C_{j_\ell}$ meets $C_{j_{\ell+1}}$ in a node for each $0\le\ell\le k_i-1$, and there are no repetitions in $j_0,\ldots,j_{k_i}$.  For $0\le\ell\le k_i-1$, let $T_\ell$ be the subset of $[n-1]$ determined by the node joining $C_{j_\ell}$ to $C_{j_{\ell+1}}$.  Then $m(g)$ acts on the restriction of $L$ to the component $C_i$ by \begin{equation}
\label{g_action_on_cpt}
g_0\prod_{\ell=0}^{k_i-1} g_{T_\ell}.
\end{equation}

Let $g\in G_x$, and suppose that $m(g)$ is trivial.  Since $m(g)$ acts on $C_0$ by $g_0$, it follows that $g_0$ is trivial.  By inducting on the number of nodes separating a given node from $C_0$, and using (\ref{g_action_on_cpt}), it follows that $g_T$ is trivial for each $T$.  Therefore, $g$ is trivial.

Now suppose $h\in G_{m(x)}$.  Let $g_0$ be the element of $\mu_r$ by which $h$ acts on $C_0$.  By inducting over the nodes as in the previous paragraph, and using the fact that the irreducible components of $C_x$ form a tree, we can now define $g_T$ for each $T$ in such a way that $m(g)=h$.  Therefore, $m:G_x\to G_{m(x)}$ is an isomorphism.\end{Prf}

We remark that \cite[section 7]{ACV} contains a careful treatment of automorphism groups for $G$-covers.

\begin{Lem}
\label{isom_criterion}
Let $\XX$ and $\YY$ be normal, separated, integral, Deligne-Mumford stacks of
finite type over an algebraically closed field $k$ of characteristic zero.
 Let $f \colon \XX\to\YY$ be a
birational morphism which induces an equivalence of categories between
objects over $Spec\;k$.  Then $f$ is an isomorphism.
\end{Lem}

\begin{Prf}
Let $V\to\YY$ be an \'etale surjective morphism from a scheme $V$ and let
$U=V\times_{\YY}\XX$.  Then $U\to V$ is separated and quasi-finite, hence
quasi-affine by \cite[A.2]{L-MB}.  Therefore, $U$ is a scheme.  Let
$U'\subseteq U$ be a connected component and let $V'\subseteq V$ be its image.
The hypotheses of the lemma imply that $U'$ and $V'$ are normal varieties and
that $U'\to V'$ is a birational morphism which is bijective on $k$-points.  By
Zariski's birational correspondence theorem, it follows that $U'\to V'$ is an
isomorphism.  Applying the argument to each connected component shows that
$U\to V$ is an isomorphism.  It now follows from \cite[3.8.1]{L-MB} that
$\XX\to\YY$ is an isomorphism.
\end{Prf}

\subsection{Comments on the construction}
The pull-back of a one-dimensional $B\mu_r$-representation is
a power of $L$. Hence, in order to understand the Chern class of the
obstruction bundle $R^1 \pi_* \C^N$ it is sufficient to understand
the Chern classes of the higher direct image $R^1 \pi_* L^w$ of powers
of $L$, and their products.

It is worth pointing out that while the ghost automorphism groups are
isomorphic to $\mu_{r_T}$, this isomorphism is not natural; the
ghost automorphism group is naturally isomorphic to the relative stabilizer
group of the twisted node, and by choosing one of the two components
$D^{T, n+1}$ or $D^{T^C, n+1}$ (and identifying the stabilizer group
by its action on the corresponding tangent bundle) one gets an
isomorphism to $\mu_{r_T}$ whose sign depends on this choice.
In our construction, this choice shows up in the definition
of $L_1$, for which we had to choose
one of $T$ and $T^C$ for all divisors $D^T$ of $\Mbar_{0, n}$. As
$L_2$ depends on that choice only up to an $r$-th power, neither
$\Mbar^{(2)}$ nor the universal line bundle $L$ depend on this choice.
The map $m \colon \Mbar^{(2)} \to \Mbar_{0, n}(e_1, \dots, e_n; B \mu_r)$
does depend on it, however. Different choices can be related by a composition
with a ghost automorphism\footnote{by which we mean an automorphism covering
the identity on the coarse moduli space} of the moduli stack.

\section{Weighted stable maps to $B\mu_r$}

\subsection{Weighted stable maps}

Let $g$ be a genus, and $\AA = (a_1, \dots, a_n)$
be \emph{weight data}, which means
$a_i \in \Q \cap [0, 1]$ satisfy
$2g - 2 + \sum_i a_i > 0$.

In \cite{Has03}, Hassett introduced the notion of weighted stable curves:
a weighted stable curve of type $(g, \AA)$ over $S$ is a nodal curve
$\pi \colon C \to S$ with $n$ sections $s_i \colon S \to C$ such that
\begin{enumerate}
\item every section $s_i$ with positive weight $a_i$ is contained in
the smooth locus of $\pi$,
\item \label{stability}
the rational divisor $K_{C/S} + \sum_i a_i s_i$ is $\pi$-relatively ample,
and
\item
for any $I \subset \otn{n}$ such that the intersection
$\bigcap_i s_i$ is non-empty, we have $\sum_i a_i \le 1$.
\end{enumerate}
We will summarize a few of his results, and refer to \cite{Has03} for
details.

If $a_i = 1$ for all $i$, then
these are stable curves in the usual sense. The difference is that when
points $s_i, i \in I$, collide, then only when $\sum_i a_i > 1$ does a new rational
component bubble off. This is enough to make the new
rational component stable according to condition (\ref{stability}).

\label{chamber-decomp}
All the moduli spaces $\Mbar_{g, \AA}$ with $|\AA|=n$ are birational. More precisely,
assume that the weight data $\AA = (a_1, \dots, a_n)$ and   $\BB = (b_1, \dots, b_n)$
satisfy $a_i \ge b_i$ for all $i$, and $a_i > b_i$ for at least one $i$ (we will
write $\AA > \BB$ from now on).  Then there is a birational reduction
morphism $\rho_{\BB, \AA} \colon \Mbar_{g, \AA} \to \Mbar_{g, \BB}$.
(It is induced by $\BB$-stabilizing the family of curves over
$\Mbar_{g, \AA}$.) There is a chamber
decomposition of $[0,1]^n$ by a finite number of walls such that the
moduli space $\Mbar_{g, \AA}$
only depends on the chamber in which the weight data $\AA$ lies:
the walls are associated to subsets $T \subset \otn{n}$ and given as
\begin{equation}			\label{eq:walls}
 w_T = \stv{a_i}{\sum\nolimits_{i \in T} a_i = 1}.
\end{equation}
Further, the contraction morphism for crossing a single wall
is given as a smooth blow-up.

It is somewhat convenient to allow at least one weight to be zero, because
$\Mbar_{g, \AA \cup \{0\}}$ is \emph{by definition} the universal curve
over $\Mbar_{g, \AA}$.

This notion has been extend to weighted stable maps in
\cite{Mustatas:Chowring}, \cite{weighted-GW} and \cite{Alexeev-Guy}.
In particular, in \cite{weighted-GW} and \cite{Alexeev-Guy}
it was shown that Gromov-Witten invariants
can be computed for any choice of weights,
yielding identical GW-invariants, and
\cite{Alexeev-Guy} gave wall-crossing formulae for the full Gromov-Witten
potential including gravitational descendants.

\subsection{$\Mbar_{0, n}$ as a blow-up of $\P^{n-3}$}
\label{sect:Pn-3-weights}
As an example that will be important later, consider for given $n$ the weights
$\AA_k = (\frac 1{k}, \dots, \frac 1{k}, 1)$ (with $n-1$ entries of
$\frac 1{k}$) for $k= 1, \dots, n-2$. 		 \label{sect:Pn-weights}
The moduli
space $\Mbar_{0, \AA_{n-2}}$ is isomorphic to $\P^{n-3}$, and the universal
curve is the blow-up $\Bl_x \P^{n-2}$
of $\P^{n-2}$ at a point $x$; the universal map
is the projection of $\P^{n-2}$ to $\P^{n-3} \cong \P T_x$ from $x$.
If we pick $x$ away from the coordinate hyperplanes, then the image
of the special section with weight 1 is the exceptional divisor, while
the remaining sections can be given as the coordinate hyperplanes.
(The special section cannot intersect with any other, while the only
condition on the remaining sections is that they may not \emph{all}
coincide.)

The moduli space $\Mbar_{0, \AA_{n-3}}$ is the blow-up of
$\Mbar_{0, \AA_{n-2}}$ at the $n-1$ points that are the images of
the intersections of $n-2$ of the $n-1$ coordinate hyperplanes.
When we successively increase the first $n-1$ weights from $\frac{1}{n-2}$
to 1, one gets a description of $\Mbar_{0, n}$ by successive blow-ups from
$\P^{n-3}$. This is also explained in \cite[section 6.2]{Has03}; the
description of $\Mbar_{0, n}$ as a blow-up of $\P^{n-3}$
is equivalent to the description by De Concini and Procesi in
\cite{deConPro:hyperplane}.

\subsection{The moduli spaces of weighted stable maps to $B\mu_r$
via $r$-th roots}
\label{subsect:weighted-stable-maps-to-Bmur}
We sidestep the question of defining a moduli problem of
weighted stable maps to a stack in general. Instead we give a direct
construction of the moduli stacks via $r$-th root constructions, guided by the
construction in the non-weighted case in section
\ref{mod-space-via-rth-roots}.

Given $r$, weight data $\AA > 0$, and
$\EE = (e_1, \dots, e_n)\in\mu_r^n$, we want to construct a
stack which would resemble $\Mbar_{0, \AA}(e_1, \dots, e_n; B \mu_r)$ if it were to exist.
Choose $d_i \in \frac 1r \Z$
with $e^{2 \pi i d_i} = e_i$ and $\sum_{i = 1}^n d_i = 0$ as before.
Boundary divisors on $\Mbar_{0, \AA}$ are given as $D_\sigma$ for
$\AA$-stable 2-partitions $\sigma = (T, T^C)$ of $\otn{n}$;
$\AA$-stable means
that the condition $\abs{T}, \abs{T^C} \ge 2$ is replaced by
$\sum_{i \in T} a_i > 1$ and $\sum_{i \in T^C} a_i > 1$. (This of
course means that a corresponding rational curve with two components
is $\AA$-stable.)

Let $\Mbar^{(1)}_{0, \AA}$ be the stack obtained from $\Mbar_{0, \AA}$
by taking the $r_T$-th root at every divisor $D^T$ such that
$(T, T^C)$ is $\AA$-stable (where $r_T$ is defined as before
as the order of $\prod_{i \in T} e_i$).
To obtain $\Cbar^{(1)}_{0, \AA}$
from $\Cbar_{0, \AA} = \Mbar_{0, \AA \cup \{0\}}$, we start with the
same construction, but additionally
construct the $r_i$-th root at every section $s_i$.\footnote{In the case of
$\Mbar_{0, n}$, the section $s_i$ is equivalent to the boundary
divisor given by $T = \{i, n+1\}$; however, this does not yield an
$\AA$-stable 2-partition, hence we need to list them separately.}
The same proof as in lemma \ref{lem:pi1} shows:
\begin{Lem}
There is a canonical map
$\Cbar^{(1)}_{0, \AA} \to \Mbar^{(1)}_{0, \AA}$.
\end{Lem}
This is a stacky curve with balanced nodes, but it can have points (in the relative
smooth locus of the coarse moduli space) with automorphism group
$\mu_r^{\abs{I}}$, for curves where $s_i, i \in I$, are identical; so this is
not a twisted stable curve in the sense of \cite{AV99}.

\begin{Prop}
Let $L_{1, \AA}$ be the line bundle on $\Cbar^{(1)}$ defined by
\begin{equation} 		\label{def:L1w}
L_{1, \AA} := \bigotimes_{i=1}^n \OO(s_i)^{d_i} \otimes
\bigotimes_T \OO(D^{T, n+1})^{d_T}
\end{equation}
where the second tensor product goes over all  subsets
$T \subset \otn{n}$ with $n \not\in T$ such that
$(T, T^C)$ is $\AA$-stable.  Then
$L_{1, \AA}^r = \pi^* L_{2, \AA}$ for some line bundle $L_{2, \AA}$ on
$\Mbar^{(1)}$.
\end{Prop}
Again, this has the same proof as before.

Let $\MbarAA$ be the $\mu_r$-gerbe over
$\Mbar^{(1)}_{0, \AA}$ of $r$-th roots of $L_{2,\AA}$,
and let $\CbarAA$ be the base change of $\Cbar^{(1)}_{0, \AA}$ to
$\MbarAA$. The line bundle
$L_\AA := L_{1, \AA} \otimes \pi^* L_{2, \AA}^{-1/r}$ on
$\overline{C}_{\AA}$ has trivial
$r$-th power and thus defines a map $f \colon \CbarAA \to B\mu_r$.

However, again $f \colon \CbarAA \to B\mu_r$ is not a twisted stable
map in the sense of \cite{AV99}; most importantly, $f$ is not representable
(not representable even after a base change to a scheme
$S \to \MbarAA$).

\subsection{$\P^{n-3}$-weight data}
\label{subsect:Pn-3-computation}

For any weight $w$ of a one-dimensional $\mu_r$-representation, we
call $$H^w_\AA = R^1 \pi_* L^w_\AA$$ the generalized dual Hodge bundle
on $\MbarAA$ for the weight data $\AA$.

Our computation of the Chern class builds up from a direct computation
for the weight data $\AA = (\frac 1{n-2}, \dots, \frac 1{n-2}, 1)$, which yields $\Mbar_{0, \AA}\cong\P^{n-3}$ (see section \ref{chamber-decomp}). Given $w$, let
$\delta_i^w \in [0,1)$ be the age of the line bundle $L^w$ at the $i$-th
section; it is determined by $e^{2 \pi i \delta^w_i} = e_i^w$. For
any subset $T \subset [n]$, we let $\delta^w_T = \sum_{i \in T} \delta^w_i$.
\begin{Prop}
The generalized dual Hodge bundle $H^w_\AA$
has the following class in the $K$-group:\footnote{See
\ref{sect:notations} for other notation conventions used in this formula.}
\begin{equation*}
\left[ H^w_\AA \right] =
\sum_{p = \langle \delta^w_{[n-1]} \rangle}^{\delta^w_{[n-1]} - 1}
\left[ \OO( - p H ) \right]
\end{equation*}
\end{Prop}

\begin{Prf}
The moduli space is a $\mu_r$-gerbe over $\P^{n-3}$, and the
universal curve is constructed from $\Bl_x \P^{n-2}$ by the ${r_i}$-th root
construction at the section $s_i$ for all $i$, and the base change to
the $\mu_r$-gerbe.

We choose $d_i$ such that $d_i \in [0,1)$
for $i = 1, \dots, n-1$ and
$d_n = - \sum_{i=1}^{n-1} d_i$. Then $L_{2, \AA}$ can be computed by
$L_2 = s_n^* (L_{1, \AA})^r = \OO(- r d_n)$, and so by projection formula
$H^w_\AA = \OO(w d_n) \otimes R^1 \pi_* L_{1, \AA}^w$. To compute the
higher direct image of $L_{1, \AA}^w$, we break up
$\pi$ into the composition $\pi = \pi_2 \circ \pi_1$ of the map
$\pi_1$, forgetting the roots along
the sections, with the map $\pi_2$ that is the base change of
the natural projection $\Bl_x \P^{n-2} \to \P^{n-3}$ to the $\mu_r$-gerbe.

The push-forward along $\pi_1$ follows easily from equation
(\ref{eq:push-forward}). Using
$-w d_n = \delta^w_{[n-1]} + \sum_{i=1}^{n-1} \floor{w d_i}$, we
get (where we write $E$ for the exceptional divisor
of $\Bl_x \P^{n-2}$):
\begin{equation*} \begin{split}
(\pi_1)_* L_{1, \AA}^w  & =
\OO\left( \sum_{i=1}^n \floor{w d_i} s_i \right)
= \OO\left( (-w d_n - \delta^w_{[n-1]})H + \floor{wd_n} E \right)
\end{split} \end{equation*}
The relative canonical bundle of $\pi_2$ is $-H - E$, and the pull-back of
the hyperplane class on $\P^{n-3}$ is $H-E$. Hence
\begin{equation*} \begin{split}
R^1 (\pi_2)_* (\pi_1)_* L_{1, \AA}^w  & =
\left((\pi_2)_* \OO\left( (w d_n + \delta^w_{[n-1]} - 1) H
			 + (-\floor{w d_n} - 1) E \right) \right)^\vee \\
 &=
\OO(- \fract{\delta^w_{[n-1]}} - wd_n) \otimes
\left((\pi_2)_* \OO\left( (\ceil{\delta^w_{[n-1]}}-2) H\right) \right)^\vee,
\end{split} \end{equation*}
where we used
$w d_n - \floor{w d_n} = \{ w d_n \}
= 1 - \fract{ - w d_n} = 1 - \fract{\delta^w_{[n-1]}}$ and
$\delta^w_{[n-1]} - \fract{\delta^w_{[n-1]}} = \ceil{\delta^w_{[n-1]}} - 1$.
Applying $(\pi_2)_*$ to
the short exact sequences $\OO( (a-1) H) \to \OO( aH) \to \OO(aH)|_H$, and using $(\pi_2)_*\OO=\OO$, implies that in the $K$-group,
$$[R^1 (\pi_2)_* (\pi_1)_* L_{1, \AA}^w] = \sum_{p=0}^{\ceil{\delta^w_{[n-1]}}-2} [\OO(-\fract{\delta^w_{[n-1]}}-wd_n-p)].$$  Tensoring this with $\OO(wd_n)$ and re-indexing yields the statement of the proposition.
\end{Prf}

Since the hyperplane class of $\P^{n-3}$ agrees with the $\psi$-class
of the $n$-th marking (which is special by having weight one), this implies:
\begin{Cor}					\label{cor:Chern-Pn}
In the situation of the previous proposition, the Chern class of
$H^2_\AA$ is given as
\[
c(H^w_\AA)
= \prod_{p = \fract{\delta^w_{[n-1]}}}^{\delta^w_{[n-1]}-1} (1 - p \psi_n)
\]
\end{Cor}
(Note that $\psi_n$ denotes the pull-back of the corresponding class
in $\Mbar_{0, \AA}$ by our convention for the rational Chow groups of
the moduli stacks.)

\subsection{Relation to the twisted $I$-function.} \label{twisted-I}
In Givental's formalism for Gromov-Witten theory
\cite{Givental:quantization, Givental:symplectic, Coates-Givental,
CCIT:computing}, the so-called
$J$-function plays an essential role. Let $X = [\C^N/\mu_r]$, where $\mu_r$ acts diagonally with weights $w_1,\ldots, w_N$.  The orbifold cohomology of $X$ is $H = H^*(X) = \bigoplus_{e \in \mu_r} \C \cdot h_e$, and the $J$ function is a map $H \to H[z][[z^{-1}]]$ defined by the following formula:
\begin{align}
J^X(z, t)
& = z + t +
\sum_{n \ge 0} \sum_{e \in \mu_r} \frac{1}{n!}
	\langle t, \dots, t, \frac{h_e}{z - \psi} \rangle^X_{0, n+1}
\cdot r h_{e^{-1}} \nonumber \\  \label{eq:J-function}
 & = z + t +
\sum_{n \ge 0} \sum_{e \in \mu_r} \sum_{k \ge 0} \frac{1}{n! z^{k+1}}
	\langle t, \dots, t, \psi_n^k h_e \rangle^X_{0, n+1}
\cdot r h_{e^{-1}}
\end{align}
For example, by the results of \cite{Jarvis-Kimura},
the $J$-function of $B \mu_r$ is given as
\[
J^{B \mu_r} (z, t) = z + t +
\sum_{\k = (k_0, k_1, \dots, k_{r-1})}
	\prod_{j = 0}^{r-1} \frac{t_j^{k_j}}{k_j! z^{\abs{k_j}}}
	\cdot h_{\prod_j (\omega^{j})^{k_j}}
\]
where we wrote $t = \sum_{j=0}^{r-1} t_j h_{\omega^{j}}$.
The idea of \cite{CCIT:computing} (and Givental's formalism in general),
specialized to our setting, is to determine the $J$-function $X$ from the
$J$-function of $B \mu_r$; the former is called the \emph{twisted
$J$-function} in \cite{CCIT:computing}.

As an approximation to the twisted $J$-function,
Coates, Corti, Iritani and Tseng define a
\emph{twisted $I$-function} in \cite[section 4]{CCIT:computing};
specialized to our case, and translated into our notation,
it is given by\footnote{The formula on page 9 ibid. defining the
``modification factor'' $M_\theta(z)$ has to be applied
with $s_i$ specialized such that $e^{\boldsymbol{s}(ch)}$ for the
Chern character $ch$ of some bundle $E$ gives the Euler class of
$-[E]$.}
\[
I^{\mathrm{tw}}(z, t) = z + t
+ \sum_{(e_1, \dots, e_n)}
J^{(e_1, \dots, e_n)} \cdot M_{(e_1, \dots, e_n)}(z). \]
Here
\[
J^{(e_1, \dots, e_n)} = z^{-n+2} \delta_{1, \prod_i e_i} \cdot h_{e_n^{-1}}
\]
is the part of the $J$-function
of $B\mu_r$ that
comes from invariants computed on $\Mbar_{0, n}(e_1, \dots, e_n; B \mu_r)$,
and $M_{(e_1, \dots, e_n)}(z)$ is defined by
\[
M_{(e_1, \dots, e_n)}(z) =
\prod_{a=1}^N
\prod_{p = \fract{\delta^{w_a}_{[n-1]}}}^{\delta^{w_a}_{[n-1]}-1} (1 - pz).
\]
To show that $I^{\mathrm{tw}}$ has the desired properties\footnote{It has
the same image as the twisted $J$-function, namely Givental's Lagrangian
cone $\LL_X$.}
they use Tseng's Grothendieck-Riemann-Roch-computation of the Chern character
of the obstruction bundle in \cite{Tseng:OQRRLS}.

We can define a \emph{weighted $J$-function} of $X$ by
\[
J^{X; \text{weighted}} (z, t) = z + t +
\sum_{n \ge 0} \sum_{e \in \mu_r} \frac{1}{n!}
	\langle t, \dots, t, \frac{h_e}{z - \psi}
	\rangle^{X; \text{weighted}}_{0, n+1}
\cdot r h_{e^{-1}},
\]
where the invariant with superscript ``weighted'' denotes the invariant
computed by the moduli of weighted stable maps
$\MbarAA$ considered in the previous section, i.e. for
$\AA = (\frac 1{n-2}, \dots, \frac 1{n-2}, 1)$. Then the
result of the previous section can be formulated as
\[
J^{X; \text{weighted}} (z, t) = z + t
+ \sum_{(e_1, \dots, e_n)}
J^{(e_1, \dots, e_n)} \cdot \widetilde M_{(e_1, \dots, e_n)}(z)
\]
where $\widetilde M_{(e_1, \dots, e_n)}(z)$ is the truncation of
$M_{(e_1, \dots, e_n)}(z)$ by $z^{n-2} = 0$:
When the Euler class of the obstruction bundle for $(e_1, \dots, e_n)$
is given as a polynomial $P(\psi_n)$ in $\psi_n$, then the contribution
to the $J$-function is the $\psi_n^{n-3}$-coefficient of
$P(\psi_n) \frac 1{z - \psi_n}$. This coefficient is given by
$z^{-n+2}$ times the truncation of $P(z)$.

Independently of $n$, this shows that $J^{X ; \text{weighted}}$ is obtained
from $I^{\mathrm{tw}}$ by removing all terms of non-negative degree in
$z$ except the first two; such terms would correspond geometrically to
a negative number of $\psi_n$-insertions.

\section{Weight change}
\label{sect:weight-change}

The goals of this section are the wall--crossing theorems
\ref{thm:K-group-wc} and \ref{thm:chern-class-wc}.

\subsection{Preparations}
We begin with several lemmata we will need in the proof.

\begin{Lem}		\label{lem:blowup-vs-rth-root}
Let $D_1, D_2$ be two smooth divisors with transversal intersection
on a smooth
Deligne-Mumford stack $X$. Let $\widetilde X$ be the blow-up of $X$ at their
intersection, with exceptional divisor $E$ and proper transforms
$\widetilde{D_1}$ and $\widetilde{D_2}$. On the other hand, consider the
$r$-th root construction $X_{(D_1, r), (D_2, r)}$ and its blowup $Z$ at the
intersection of the two gerbes $\DD_1$ and $\DD_2$ lying over $D_1$ and $D_2$,
respectively.  Then $Z$ is isomorphic to the $r$-th root construction
$\widetilde X_{(\widetilde{D_1}, r), (\widetilde{D_2}, r), (E, r)}$.
\end{Lem}

\begin{Prf}
Let $\widetilde{\DD_i}$ and $\EE$ be the gerbes in
$\widetilde X_{(\widetilde{D_1}, r), (\widetilde{D_2}, r), (E, r)}$
lying over $\widetilde{D_i}$ and $E$, respectively.
The line bundles $\OO(\widetilde{\DD}_i \otimes \EE)$
with their
canonical sections are $r$-th roots of the pull-backs of $D_i \subset X$,
determining a morphism
\[ f \colon \widetilde X_{(\widetilde{D_1}, r), (\widetilde{D_2}, r), (E, r)}
\to X_{(D_1, r), (D_2, r)}. \]
The pull-back of the ideal sheaf of
the intersection $\DD_1 \cap \DD_2$
along $f$
is the ideal sheaf of $\EE$: this is easy to see locally,
where we can assume that $\DD_i$ is cut out by an equation
$(\ss_i = 0)$; its pull-back $f^* \ss_i$ cuts out
$\EE \cup \widetilde{\DD_i}$. Since $\EE$ is Cartier,
the universal property of blow-ups yields a map
\[
g \colon \widetilde X_{(\widetilde{D_1}, r), (\widetilde{D_2}, r), (E, r)}
\to Z.
\]

To go the other way, we first show that $Z\to X$ lifts to $\tilde{X}$.  The preimage of $D_1\cap D_2$ in $X_{(D_1, r), (D_2, r)}$ is $\DD_1^r\cap\DD_2^r$.  One can check after \'etale base change to a scheme that the preimage of this in $Z$ is $r$ times the exceptional divisor, hence is Cartier.  So the universal property of blowups gives us a morphism $Z\to \widetilde X$.  The preimage of $E$ under this morphism is $r$ times the exceptional divisor of $Z$, and it follows that the preimage of $\widetilde{D_i}$ is $r$ times the proper transform of $\DD_i$.  This gives us a lifting to $$h:Z\to \widetilde X_{(\widetilde{D_1}, r), (\widetilde{D_2}, r), (E, r)}.$$

As neither $Z$ nor $\widetilde X_{(\widetilde{D_1}, r), (\widetilde{D_2}, r), (E, r)}$ have nontrivial automorphisms over the identity of $X$, both $gh$ and $hg$ must be (2-isomorphic to) the identity.

\end{Prf}

\begin{Lem} \label{lem:higher-images}
Let $X, Y$ be Deligne-Mumford stacks, and let $f \colon X \to Y$ be
a composition of $r$-th root constructions
and blow-ups at regularly embedded centers (i.e. the normal sheaf of the
center is a vector bundle). Then $\R f_* \L f^* \FF = \FF$
for any quasi-coherent sheaf $\FF$ on $Y$.
\end{Lem}
\begin{Prf}
By the projection formula, it is enough to prove $\R f_* \OO_X = \OO_Y$. For
$r$-th root constructions, it is obvious that the higher direct images
vanish, and
by \cite{chuck:using-stacks-to-impose} $f_* \OO_X = \OO_Y$. For
blow-ups, this is well-known in the case of schemes.  Since blow-ups are representable, one can reduce to the case of schemes by taking an \'etale base change to a scheme covering $Y$.
\end{Prf}

\begin{Lem}  \label{lem:base-change} (Base change)
Let
\[\xymatrix{
Y' \ar[r]^h \ar[d]_{f'} \ar@{}[dr]|{\Box} & Y \ar[d]^f \\
X' \ar[r]_g & X}\]
be a 2-cartesian square of Deligne-Mumford stacks where $X$ is
quasi-compact, $f$ is quasi-compact and quasi-separated, and $Y$ and $X'$ are tor-independent over $X$.  Let $\EE\in D^b(Y)$ be a complex having quasi-coherent cohomology.  Suppose either that $g$ has finite tor-dimension or that $\EE$ has flat, finite amplitude relative to $f$.  Then there is a natural isomorphism
$$\L g^*\R f_*\EE \to \R f_*' \L h^*\EE.$$
\end{Lem}
We will apply this in the case where $f$ is flat and $\EE$ is a line
bundle on $Y$.

\begin{Prf}
The existence of a natural morphism $\L g^*\R f_*\EE \to \R f_*' \L h^*\EE$
follows from the adjointness of pullback and pushforward.  Indeed, the natural
morphism $\EE \to \R h_*\L h^*\EE$ determines a morphism $$\R f_*\EE \to \R
(fh)_*\L h^* \EE = \R (gf')_* \L h^*\EE,$$ which is equivalent to the
morphism above by adjointness of $\R g_*$ and $\L g^*$.

For schemes, the proposition is the same as \cite[IV, 3.1.0]{SGA6}, and the
reduction to the case of schemes is identical to the proof of
\cite[13.1.9]{L-MB}.
\end{Prf}

\begin{Lem} 				\label{lem:blowup-pushforward}
Let $\tau \colon \widetilde X \to X$ be the blow-up of a smooth
Deligne-Mumford stack $X$ at a smooth center $Z \subset X$
of codimension two, with normal bundle $N$ and exceptional divisor $E$.
Then for $n \ge 0$
\[
[R^1 \tau_* \OO(n \cdot E)]
= \sum_{k=0}^{n-2} [\Lambda^2 N \otimes \Sym^k N].
\]
\end{Lem}
\begin{Prf}
As a blow-up is representable, it is sufficient to check this for schemes.
Due to the short exact sequences
$\OO((n-1) \cdot E) \to \OO(n E) \to \OO_E(-n)$ on $\widetilde X$,
it follows from induction
if we show $R^1 \pi_* \OO_E(-n) = [\Lambda^2 N \otimes \Sym^{n-2} N]$.
This is easily checked by Serre duality, since $$R^0\pi_*\OO_E(n) = \Sym^n(N^\vee)$$ and the relative dualizing sheaf of $\tau|_E$ is $\OO_E(-2)\otimes\tau^*(\Lambda^2N)^{-1}$.
\end{Prf}

\subsection{Constructing the reduction map}

Consider two weight data $\AA > \BB$; then there is a reduction
morphism $\rho_{\BB, \AA} \colon \Mbar_{0, \AA} \to \Mbar_{0, \BB}$.
Further, we assume the following property:
\begin{itemize}
\item[(*)]
There is exactly one 2-partition $\sigma = (T_0, T_0^C)$ of $\otn{n}$
such that $\sigma$ is $\AA$-stable but $\BB$-unstable.
\end{itemize}
In other words, there is just one
wall between $\AA$ and $\BB$ in the chamber decomposition of the set
of weight data discussed in section \ref{chamber-decomp}.
Specifically, this means that
$\sum_{i \in T_0} a_i > 1$ but $\sum_{i \in T_0} b_i \le 1$, and that $T_0$
is the only such subset of $\otn{n}$.

Then $\rho_{\BB, \AA} \colon \Mbar_{0, \AA} \to \Mbar_{0, \BB}$ is the
blow-up of $\Mbar_{0, \BB}$ at the locus $Z_{T_0}$
of curves where all $s_i, i \in T_0$ agree; the
exceptional divisor is $D^{T_0}$. The universal curve
$\Cbar_{\AA}$ is obtained from $\Cbar_{\BB}$ in two steps:
\begin{enumerate}
\item Blowing up
at the preimage $\pi^{-1}Z_{T_0}$ of $Z_{T_0}$, i.e. taking the base change
for $\rho_{\BB, \AA}$; we denote the resulting family over
$\Mbar_{0, \AA}$ by $\Cbar_{\BB/\AA}$.
\item Blowing up at the preimage
(with respect to the previous blow-up) of
the common image $s_i(Z_{T_0})$ for any $i \in {T_0}$. (See \cite[Remark
3.1.2]{weighted-GW}.)
\end{enumerate}

There is a canonical map
\[ \rho_{\BB, \AA}^{(1)} \colon
\Mbar^{(1)}_{0, \AA} \to \Mbar^{(1)}_{0, \BB}. \]
This follows from the fact that the pull-back of any boundary divisor
$D^S \subset \Mbar_{0, \BB}$ is just the
corresponding divisor $D^S \subset \Mbar_{0, \AA}$, and that the $r$-th root
construction commutes with such base change.

A matching map
$\Cbar^{(1)}_{\AA} \to \Cbar^{(1)}_{\BB}$ does not
exist in general, for the following reason: the preimage of
a section $s_i \subset \Cbar_{\BB}$ with $i \in {T_0}$ is the union
of the corresponding section $s_i \in \Cbar_{\AA}$ with the exceptional
divisor $D^{T_0, n+1}$ of the second blow-up step in the construction
of $\Cbar_{0, \AA}$ above.  However, for example
when $\prod_{i \in T_0} e_i = 1$, there is no root construction
along the divisor $D^{T_0, n+1}$ in the construction of $\Cbar^{(1)}_{\AA}$
at all, and so an $r_i$-th root of the pull-back of
$s_i \subset \Cbar_{0, \BB}$ does not exist.

In order to compare the Hodge bundles,
we will later construct some auxiliary
spaces to overcome this problem.

For simplicity, we make the following additional assumption:
\begin{itemize}
\item[(**)]
We assume that for ${T_0}$ as in (*), we have $n \not\in {T_0}$.
\end{itemize}
(This simplifies the computation with respect to our choice of
$L_{1, \AA}$ in definition \ref{def:L1w} and, on the other hand,
always holds when we start with the $\P^{n-3}$-weight data used in
section \ref{subsect:Pn-3-computation}.)

\begin{Lem}
Assuming (*) and (**), we have $\rho_{\BB, \AA}^* L_{2, \BB} = L_{2, \AA}$
\end{Lem}
This is immediate from $L_{2, \AA} = s_n^* (L_{1, \AA})^r$, as $s_n$ does
not meet any of the divisors appearing in the definition (\ref{def:L1w}) of
$L_{1, \AA}$, except itself.  As a consequence, it follows that:
\begin{Cor}
There is a well-defined reduction map
\[
\rho_{\BB, \AA} \colon \MbarAA \to \MbarBB.
\]
\end{Cor}

\subsection{Weight change and Hodge bundles}

We continue with the assumptions (*) and (**) from the previous section.
We want to compare $\rho_{\BB, \AA}^* H^w_\BB$ and $H^w_\AA$ in the
$K$-group of $\MbarAA$, where
due to projection formula $H^w_\AA$ can be computed by
\[ H^w_\AA = R^1 \pi_* L_{1, \AA}^w \otimes \left(L_{2, \AA}\right)^{-\frac wr}.
\]

We will introduce several auxiliary spaces; the goal is to have a
reduction map as a smooth blow-up between spaces that are very close
to the universal curves. This is achieved in the map $\tau$ below.

Let $A_{\AA}$ be the intermediate space obtained from $\CbarAA$
by forgetting the root construction along all sections. In other words,
it is constructed from the universal curve $\Cbar_\AA$ by the $r_T$-th root
construction for every divisor $D^T$ where $T \subset \otn{n+1}$ and
$T$ is $\AA \cup \{0\}$-stable, followed by the base change along
$\MbarAA \to \Mbar_{0, \AA}^{(1)}$.
Let $\nu_\AA \colon \CbarAA \to A_\AA$ be the induced map,
and $\pi_\AA' \colon A_\AA \to \MbarAA$ the projection
to the moduli space.

Then
\begin{equation} \label{nu*}
 {\nu_\AA}_{*} \left(L_{1, \AA}^w\right)
= \OO \left( \sum_{i=1}^n \floor{w d_i} s_i +
\sum_{\substack{T \subset [n-1] \\
\text{$(T, T^C)$ is $\AA$-stable}}} w d_T \cdot D^{T, n+1} \right)
\end{equation}

Now let $A_{\BB/\AA}$ be the base change of
$A_\BB \to \MbarBB$ along $\rho_{\BB, \AA}$. While
a map $A_\AA \to A_{\BB/\AA}$ exists, we prefer not to use it and
instead consider two more additional spaces:
Pick any section $s_{j_0}$ with ${j_0} \in T_0$ and let
$A'_\AA = \left(A_\AA\right)_{s_{j_0}, r_T}$ be the stack obtained from
$A_\AA$ by adding the $r_T$-th root construction at the ${j_0}$-th section.
We define $A'_{\BB/\AA}$ analogously.

\[
\xymatrix{
\CbarAA \ar[rd]^{\nu_\AA} &
{A'_\AA} \ar[d]^{\epsilon_\AA} \ar[r]^\tau
& A'_{\BB/\AA} \ar[d]^{\epsilon_{\BB/\AA}}
& & \CbarBB \ar[ld]^{\nu_\BB} \\
& {A_\AA} \ar[rd]_{\pi'_\AA}
& {A_{\BB/\AA}} \ar[r]^{\rho_{\BB, \AA}} \ar[d]^{\pi'_{\BB/\AA}}
& {A_\BB} \ar[d]^{\pi'_\BB} \\
& & {\MbarAA} \ar[r]^{\rho_{\BB, \AA}} & {\MbarBB}
}
\]

Applying lemma \ref{lem:blowup-vs-rth-root} to the divisors
$s_{j_0}$ and $\pi^{-1} D^T$ on the coarse moduli space of $A_{\BB/\AA}$,
we see that the map $\tau$ is the blow-up
at the intersection of the tautological gerbes over $s_{j_0}$ and
$\left(\pi'_{\BB/\AA}\right)^{-1} D^{T_0}$ in $A'_{\BB/\AA}$.

Setting $L^w_{1, \BB/\AA} = \rho_{\BB, \AA}^* {\nu_\BB}_* L^w_{1, \BB}$, we have
\begin{align}
(\rho_{\BB, \AA})^* H^w_\BB =&
(\rho_{\BB, \AA})^* R^1 (\pi'_{\BB})_* {\nu_\BB}_* L_{1, \BB}^w
\otimes L_{2, \BB}^{-\frac wr} & \nonumber \\
=&  R^1 (\pi'_{\BB/\AA})_* \rho_{\BB, \AA}^* {\nu_\BB}_* L_{1, \BB}^w
\otimes (L_{2, \AA})^{-\frac wr} &
\text{(lemma \ref{lem:base-change})} \nonumber \\
=&  R^1 (\pi'_{\BB/\AA})_* L_{1, \BB/\AA}^w
\otimes (L_{2, \AA})^{-\frac wr} &
\label{H^w_B} \\
=&  R^1 (\pi'_{\BB/\AA} \epsilon_{\BB/\AA} \tau)_*
(\epsilon_{\BB/\AA} \tau)^* L_{1, \BB/\AA}^w
\otimes (L_{2, \AA})^{-\frac wr} &
\text{(lemma \ref{lem:higher-images})} \nonumber \\
\intertext{on the other hand,}
H^w_\AA =&
R^1 (\pi'_{\AA})_* {\nu_\AA}_* L_{1, \AA}^w \otimes (L_{2, \AA})^{-\frac wr} \nonumber \\
=& R^1 (\pi'_\AA \epsilon_\AA)_* \epsilon_\AA^* {\nu_\AA}_* L_{1, \AA}^w
\otimes (L_{2, \AA})^{-\frac wr}
& \text{(lemma \ref{lem:higher-images})}\nonumber \\
=& R^1 (\pi'_{\BB/\AA} \epsilon_{\BB/\AA} \tau)_*
\epsilon_{\AA}^* {\nu_\AA}_* L_{1, \AA}^w
\otimes (L_{2, \AA})^{-\frac wr}. \label{H^w_A}
\end{align}

Now for any $T \subset [n-1]$ such that $(T, T^C)$ is
$\BB$-stable, it is also $\AA$-stable, and
\[ (\rho_{\BB, \AA} \epsilon_{\BB/\AA} \tau)^* (D^{T, n+1})
= \epsilon_\AA^* (D^{T, n+1}).
\]
The pull-back of the divisor class of sections is given by
\[
(\rho_{\BB, \AA} \epsilon_{\BB/\AA} \tau)^* (s_i)
= \begin{cases} s_i & \text{if $i \not \in T_0$} \\
s_i + D^{T_0, n+1} & \text{if $i \in T_0$}
\end{cases}
\]
Using these formulae and equation (\ref{nu*}) for $\AA$ and $\BB$,
respectively, yields
\[
(\epsilon_\AA)^* {\nu_\AA}_* L_{1, \AA}^w
= (\epsilon_{\BB/\AA} \tau)^* L_{1, \BB/\AA}^w \otimes
\OO \left( \sum_{i \in T_0} \bigl(wd_i - \floor{wd_i}\bigr)\cdot D^{T_0, n+1}
\right). \]
Note that $\delta^w_{T_0} = \sum_{i \in T_0} \bigl(wd_i - \floor{wd_i}\bigr)$.
The projection formula yields
\begin{equation}\label{wall_crossing_step}
\R (\epsilon_{\BB/\AA} \tau)_* (\epsilon_\AA)^* {\nu_\AA}_* L_{1, \AA}^w
= L_{1, \BB/\AA}^w \otimes
\R (\epsilon_{\BB/\AA} \tau)_* \OO \bigl(\delta^w_{T_0} \cdot
D^{T_0, n+1}\bigr).
\end{equation}
Combining equations (\ref{H^w_B}, \ref{H^w_A}, \ref{wall_crossing_step}), and using
$ R^0 (\epsilon_{\BB/\AA} \tau)_* \OO \bigl(\delta^w_{T_0} \cdot
D^{T_0, n+1}\bigr)
= \OO$, it follows that
\begin{equation}\label{wall_crossing_step2}
[H^w_\AA] = [\rho_{\BB, \AA}^* H^w_\BB] +
\left[(\pi'_{\BB/\AA})_* \left( R^1(\epsilon_{\BB/\AA} \tau)_*
\OO \bigl(\delta^w_{T_0} D^{T_0, n+1}\bigr) \otimes L_{1, \BB/\AA}^w \right)
\otimes L_{2, \AA}^{- \frac wr}\right]
\end{equation}
Write $\delta^w_{T_0}$ as the fraction $\frac {p_0}{r_0}$, where $r_0 =
r_{T_0}$.
The exceptional divisor of $\tau$ is
$\DD^{T_0, n+1} = \frac 1{r_0} D^{T_0, n+1}$; so
by lemma \ref{lem:blowup-pushforward}
\begin{equation} 		\label{eq:discrepancy-pushforward}
[R^1 \tau_* \OO\bigl(\delta^w_{T_0} \cdot D^{T_0, n+1}\bigr)]
= \sum_{k=0}^{p_0 - 2} [\Lambda^2 N \otimes \Sym^k N],
\end{equation}
where $N$ is the normal bundle to the center of the blow-up $\tau$.
To compute the right-hand side, we introduce additional
normal bundles. Let $N_{\ss_{j_0}}$ be the normal bundle of the gerbe
$\ss_{j_0}$ over $s_{j_0}$ in $A'_{\BB/\AA}$, let
$N_{\DD^{T_0}}$ be the normal bundle to the gerbe $\DD^{T_0}$
over $D^{T_0}$ in $\MbarAA$, and let $N_{s_{j_0}}$ be the normal bundle
to $s_{j_0}$ in $A_{\BB/\AA}$ (equivalently, $N_{s_{j_0}}$ is the relative
tangent bundle of $\pi'_{\BB/\AA}$ restricted to $s_{j_0}$).
Then $N_{\ss_{j_0}}$ is the restriction of the tautological bundle of the
$r_0$-th root construction
$\epsilon_{\BB/\AA}$ at the section $s_{j_0}$ to the gerbe, and
\begin{equation} 			\label{eq:eps-pushforward}
 (\epsilon_{\BB/\AA})_* N_{\ss_{j_0}}^k =
\begin{cases}
N_{s_{j_0}}^{\frac k{r_0}} & \text{if $r_0$ divides $k$,} \\
0 & \text{otherwise.}
\end{cases}
\end{equation}
The section $s_{j_0}$ induces a splitting of the tangent bundle
of $A_{\BB/\AA}$ along $s_{j_0}$ into the relative tangent bundle
and the push-forward of the tangent bundle of $\MbarAA$ along $s_{j_0}$.
This induces a splitting of $N$ as
\[
N = \restr{\left((\epsilon_{\BB/\AA} \pi'_{\BB/\AA})^* N_{\DD^{T_0}}
	\oplus N_{\ss_{j_0}}\right)}{Z_{T_0}}
\]
Applying this splitting to equation (\ref{eq:discrepancy-pushforward}),
we obtain
\[
[R^1 \tau_* \OO\bigl(\delta^w_{T_0} \cdot D^{T_0, n+1}\bigr)]
= \sum_{k=0}^{p_0 - 2} \sum_{a+b = k}
[N_{\ss_{j_0}}^{b+1} \otimes (\epsilon_{\BB/\AA} \pi'_{\BB/\AA})^*
N_{\DD^{T_0}}^{a+1} ].
\]
We can rewrite the summation as $\sum_{b=0}^{p_0-2}\sum_{a=0}^{p_0-b-2}$ or,
equivalently, $\sum_{b=0}^{p_0-1}\sum_{a=0}^{p_0-b-2}$ 
By equation (\ref{eq:eps-pushforward}), only the terms with
$m r_0 =b+1$ for some $m \in \Z$ are surviving the push-forward
along $\epsilon_{\BB/\AA}$, which yields:
\[
[R^1 (\epsilon_{\BB/\AA} \tau)_*
	\OO\bigl(\delta^w_{T_0} \cdot D^{T_0, n+1}\bigr)]
= \sum_{m=1}^{\ceil{\delta^w_{T_0}}-1} \sum_{a = 0}^{p_0 - r_0m-1}
[N_{s_{j_0}}^m  \otimes (\pi'_{\BB/\AA})^* N_{\DD^{T_0}}^{a+1}]
\]
Combining this with equation (\ref{wall_crossing_step2}) yields the wall--crossing theorem in the $K$-group:
\begin{Thm} 				\label{thm:K-group-wc}
The generalized dual Hodge bundles
$H^w_\AA$ and $H^w_\BB$ can be related in the $K$-group of $\MbarAA$ as
follows:
\begin{equation} 		\label{eq:K-group-wc}
[H^w_\AA] = [\rho_{\BB, \AA}^* H^w_\BB] +
\sum_{m=1}^{\ceil{\delta^w_{T_0}}-1} \sum_{a = 1}^{p_0 - r_0m}
[N_{\DD^{T_0}}^a \otimes
s_{j_0}^* N_{s_{j_0}}^m  \otimes s_{j_0}^* L_{1, \BB/\AA}^w
\otimes L_{2, \AA}^{-\frac wr}]
\end{equation}
where $s_{j_0}$ is the section $s_{j_0} \colon \MbarAA \to A_{\BB/\AA}$.
\end{Thm}

Let $\alpha$ be the first Chern class of the line bundle
$s_{j_0}^* {\nu_\BB}_* L_{1, \BB}^w \otimes L_{2, \BB}^{- \frac wr}$ on $\MbarBB$.
The first Chern class of $s_{j_0}^* N_{s_{j_0}}$ on $\MbarBB$
is $- \psi_{j_0}$, and so the first
Chern class of
$\KK := s_{j_0}^* N_{s_{j_0}}^m \otimes s_{j_0}^* L_{1, \BB/\AA}^w
	\otimes L_{2, \AA}^{-\frac 2r}$
is $\rho_{\BB, \AA}^* (\alpha - m \psi_{j_0})$.

Tensoring the short exact sequences
\[
0 \to \OO(\frac 1{r_0} D^{T_0})^{a-1}
  \to \OO(\frac 1{r_0} D^{T_0})^a
  \to N_{\DD^{T_0}}^a \to 0
\]
with $\KK$ yields Chern classes for all summands of the
right hand side of equation (\ref{eq:K-group-wc}), and thus the following
formula relating the Chern classes of the dual Hodge bundles:
\begin{align*}
c(H^w_\AA) & = \rho_{\BB, \AA}^*(c(H^w_\BB)) \cdot
\prod_{m=1}^{\ceil{\delta^w_{T_0}}-1} \prod_{a = 1}^{p_0 - r_0m}
\frac {1 + \frac a{r_0} D^{T_0} + \rho_{\BB,\AA}^* (\alpha - m \psi_{j_0})}
      {1 + \frac {a-1}{r_0} D^{T_0} + \rho_{\BB,\AA}^* (\alpha - m
\psi_{j_0})} \\
& = \rho_{\BB, \AA}^*(c(H^w_\BB)) \cdot
\prod_{m=1}^{\ceil{\delta^w_{T_0}}-1}
\frac {1 + (\delta^w_{T_0} - m) D^{T_0} + \rho_{\BB,\AA}^* (\alpha - m
\psi_{j_0})}
      {1 + \rho_{\BB,\AA}^* (\alpha - m \psi_{j_0})} \\
& = \rho_{\BB, \AA}^*(c(H^w_\BB)) \cdot
\prod_{m=1}^{\ceil{\delta^w_{T_0}}-1}
\left(1+\frac{(\delta^w_{T_0} - m) D^{T_0}}
{1 + \rho_{\BB,\AA}^* (\alpha - m \psi_{j_0})}\right)
\end{align*}

Equation (\ref{eq:K-group-wc}) implies that this formula does not depend on $\alpha$ itself but only on
its restriction $\alpha \cdot Z_{T_0}$ to the center of the blowup-part
of $\rho_{\BB, \AA}$.
The Chern class of $L_{1, \BB}$ on $\CbarBB$ is by
construction equal to the Chern class of $\pi_\BB^* L_{2, \BB}^{\frac 1r}$.
Hence $\alpha$  can be computed as the difference of the Chern classes
of ${\nu_\BB}_* L_{1, \BB}^w$ and $L_{1, \BB}^w$, which is given by
$\sum_{i = 1}^n (\floor{wd_i} - wd_i) s_i = -\sum_{i=1}^n \delta^w_i s_i$.
Pulling this back via $s_{j_0}$, and ignoring everything that
restricts as zero to $Z_{T_0}$, gives
$\alpha = -\sum_{i \in T_0} \delta^w_i (- \psi_{j_0})
= \delta^w_{T_0} \psi_{j_0}$.

We claim that the restriction of
$\rho_{\BB, \AA}^* \psi_{j_0}$ to $D^{T_0}$ is the $\psi$-class of the
node on the component corresponding to the complement of $T_0$. To see this, note that the irreducible component over $D^{T_0}$ corresponding to $T_0^C$ is the pull-back of the family over $Z_{T_0}$, with $s_{j_0}$ being pulled back to the node.

Since in the final formula, after expanding the fraction using a geometric
series, $\rho_{\BB, \AA}^* \psi_{j_0}$ only appears in monomials that also
have a factor of $D^{T_0}$, we can replace $\rho_{\BB, \AA}^* \psi_{j_0}$
by $\psi_{T_0}$; here $\psi_{T_0}$ is for now any divisor on that restricts
as the $\psi$-class of the node to $D^{T_0}$ (but see (\ref{eq:psi_T}) for
a somewhat canonical global definition).
So the formula simplifies further:
\begin{Thm}			\label{thm:chern-class-wc}
Assume that there is a single wall $w_{T_0}$ as defined in
(\ref{eq:walls}) between the two weight data $\AA \ge \BB$. The Chern classes
of the generalized dual Hodge bundles can be related on $\MbarAA$ as follows:
\[
c(H^w_\AA)  = \rho_{\BB, \AA}^*(c(H^w_\BB)) \cdot
\prod_{p=\fract{\delta^w_{T_0}}}^{\delta^w_{T_0} - 1}
\left(1 +  \frac {p D^{T_0}} {1 + p \psi_{T_0}}\right)
\]
\end{Thm}

\section{Chern class formula}

\subsection{Main theorem}
Theorem \ref{thm:chern-class-wc} and corollary \ref{cor:Chern-Pn}
immediately give a closed formula for the equivariant Euler class
of the generalized dual Hodge bundle.

Assume that $\mu_r$ is acting diagonally on $\C^N$ with weights
$w_1, \dots, w_N$.  Given $e_1, \dots, e_n \in \mu_r$,
let $\delta_i^{(a)} \in [0,1)$
be the age of $e_i$ acting on the $a$-th coordinate direction, i.e.
$e^{2 \pi i \delta_i^{(a)}} = e_i^{w_a}$.
For
all subsets $T \subset [n]$, let
$\delta_T^{(a)} = \sum_{i \in T} \delta_i^{(a)}$.
\begin{Thm} 		\label{mainthm}
On the connected component $\Mbar_{0, n}(e_1, \dots, e_n; B \mu_r)$
of the moduli space of twisted stable maps $\Mbar_{0, n}(B\mu_r)$,
the equivariant Euler class of the obstruction
bundle is given as
\[
e_\TT \left([R^1 \pi_* f^* \C^N]\right)
= \prod_{a=1}^N
\prod_{p= \fract{\delta^{(a)}_{[n-1]}}}^{\delta^{(a)}_{[n-1]} - 1}
(t_a - p \psi_n) \cdot
\prod_{\substack{T \subsetneq \otn{n-1} \\ 2 \le \abs{T}}}
\prod_{p = \fract{\delta^{(a)}_T}}^{\delta^{(a)}_T - 1}
	\left(1 + \frac{ p D^T}{t_a + p \psi_T} \right)
\]
\end{Thm}

\begin{Prf}
If we start with weight data $\AA = (\frac 1{n-2}, \dots, \frac 1{n-2}, 1)$
as in section \ref{subsect:Pn-3-computation}, we can choose a path in
$[0,1]^n$ leading to $\AA = (1, 1, \dots, 1)$ such that we pass every wall
$w_T = \stv{a_i}{\sum_{i \in T} a_i = 1}$ for
$T \subsetneq [n-1], \abs{T} \ge 2$ exactly once, and only one wall
at a time. By theorem \ref{thm:chern-class-wc},
we pick up exactly the factor in the above product corresponding to $D^T$
when we cross the wall $w_T$, after we set $t_a = 1$.
To get the equivariant Euler class from the total Chern class, we just have
to multiply the $i$-th Chern class of the higher direct image of the
$a$-th coordinate direction $\C \subset \C^N$
with $t_a^{\rk - i}$, as the torus is acting trivially
on the moduli space, and linearly with multiplication by $t_a$ on the
fibers of the vector bundle.
\end{Prf}

\subsection{Remarks on the formula}
				\label{subsect:formula-remarks}
Using the notations $D^{[n-1]} = - \psi_n$, $D^{i} = -\psi_i$ and
$\psi_{[n-1]} = 0$ as explained in the appendix, the formula can be written in
more compact form:
\begin{equation}		\label{eq:chern-class-compact}
e_\TT \left([R^1 \pi_* f^* \C^N] - [R^0 \pi_* f^* \C^N] \right)
= \prod_{a=1}^N t_a^{\delta_{[n]}^{(a)}-1}
\prod_{\emptyset \neq T \subseteq \otn{n-1}}
\prod_{p = \fract{\delta^{(a)}_T}}^{\delta^{(a)}_T - 1}
	\frac{t_a + p \psi_T + p D^T}{t_a + p \psi_T}
\end{equation}
Here $\delta_{[n]}^{(a)}-1$ is the virtual dimension of the contribution
from the $a$-th
coordinate direction $\C \subset \C^N$ to obstruction bundle.

The convenience of this formulation is that it remains correct
(up to an overall power of $t_a$) as long as
$\delta^{(a)}_{\{i\}}$ is any real number such that
$e_i^{w_a} = e^{2 \pi i \delta^{(a)}_{\{i\}}}$
(if we still define $\delta_T^{(a)} = \sum_{i \in T}
\delta_{\{i\}}^{(a)}$). This is shown in the appendix, see
lemma \ref{combinatorial_relation}.

This version of the formula also gives the correct answer for the
necessary localization computation in case all $e_i$ act trivially
on one of the coordinate direction, that is if $e_i^{w_a} = 1$ for
some $a$ and all $i$. In that
case, $\Mbar_{0, n}(e_1, \dots, e_n; B \mu_r)$ is the fixed
point locus of $\Mbar_{0, n}(e_1, \dots, e_n; [\C^N/\mu_r])$ (instead
of being isomorphic to it). The factor of $\frac 1{t_a}$ we get in the above
formula is the contribution of the $a$-th coordinate direction to the inverse
of the equivariant Euler class of the normal bundle of the fixed point locus.

\section{Recursions for Gromov-Witten invariants}
\label{sect:recursions}

\subsection{Inclusion-exclusion principle}
The formula gives particularly nice recursions when the invariants are
(almost) non-equivariant. To expand the formula, we use the following
fact, which we think of as a generalized inclusion-exclusion principle:
\begin{Lem}				\label{lem:inclusion-exclusion}
Let $S$ be a partially ordered set. Let $U(S)$ be the set of non-empty
subsets $I\subset S$ such that no two elements of $I$ are comparable. For
every subset $I \subset S$, let $C(I) \subset S$ be the
``ordered complement'' of $I$: the set of elements of $S$ that are not
less than or equal to any element of $I$. Then:
\begin{equation} 		\label{eq:inclusion-exclusion}
\prod_{T \in S} (1 + x_T) = 1 +
\sum_{I \in U(S)} (-1)^{\abs{I} + 1} \prod_{T \in I} x_T
\prod_{T \in C(I)} (1 + x_T)
\end{equation}
\end{Lem}
\begin{Prf}
For any subset $J \subset S$, the monomial $\prod_{T \in J} x_T$
appears in the right-hand side of the above product whenever $I$ is a
subset of the set of minimal elements of $J$. It is easily checked that
it overall has coefficient one.
\end{Prf}

We use this for $S$ being the set of subsets $T \subset [n-1]$ with
$2 \le \abs T \le n-2$, ordered by inclusion, and
\begin{align*}
x_T &= -1 + \prod_{a=1}^N
\prod_{p = \fract{\delta^{(a)}_T}}^{\delta^{(a)}_T - 1}
	\(1 + \frac{ p D^T}{t_a + p \psi_T} \) \\
&= \frac{\prod_{a,p} (1-t_a^{-1}p\psi_{T^C}) - \prod_{a,p} (1+t_a^{-1}p\psi_T)}{\prod_{a,p} (1+t_a^{-1}p\psi_T)},
\end{align*}
where $\psi_{T^C}=-D^T-\psi_T$ (see Appendix).
Then $x_T$ is a class with support on $D^T$. Thus, if $T, T'$ are not
comparable, then $x_T x_{T'}$ can only be non-zero if $T$ and $T'$ are
disjoint; hence the expansion of lemma \ref{lem:inclusion-exclusion} reduces
to a sum over a combination of pairwise disjoint subsets
$T_1, \dots, T_k \subset [n-1]$. Since every $T_i$ has size at least
2, the datum of $\{T_1, \dots, T_k\}$ can
be identified with a partition $\PP$ of $[n-1]$; given $\PP$,
we write $\PP_{\ge 2}$ for the sets
in $\PP$ that have size at least 2, recovering the list of the $T_i$.
We can similarly simplify the second product of equation
(\ref{eq:inclusion-exclusion}) to a product over $T$ which are either disjoint from
or fully contain $T_i$, for all $i$; in other words, we can identify
$T$ with a subset of the quotient set
$[n-1]/(\PP)$ having at least $2$ elements. Thus:

\begin{multline} 			\label{eq:incl-excl-applied}
e_\TT (R^1 \pi_* f^* \C^N) = \\
\left(
	1 + \sum_{\PP} (-1)^{\abs{\PP_{\ge 2}}+1}
	\prod_{T \in \PP_{\ge 2}} x_T
	\prod_{\substack{T \subsetneq [n-1]/(\PP) \\ 2 \le \abs{T}}}
		(1 + x_T)
\right)
\prod_{a=1}^N
\prod_{p= \fract{\delta^{(a)}_{[n-1]}}}^{\delta^{(a)}_{[n-1]} - 1}
(t_a - p \psi_n)
\end{multline}
where the sum goes over all non-trivial partitions
$\PP$ of $[n-1]$ (excluding the partitions of size 1 and $n-1$), and
we identify a subset $T\subset [n-1]/(\PP)$ with its preimage in $[n-1]$.

The class associated to the partition $\PP$ in the above expansion has
support on $\bigcap_{T \in \PP_{\ge 2}} D^T$,
which explains why we call it an inclusion-exclusion principle.
This intersection is a moduli space of comb curves as in the figure
on page \pageref{fig:recursion}.

\subsection{Non-equivariant recursions for $[\C^3/\mu_r]$}\label{s:non_eq_c3mur}
Let $\mu_r$ act non-trivially on $\C^3$ so that it leaves the volume form
of $\C^3$ invariant. Up to isomorphism of $\mu_r$, we may assume that the
generator is acting with age 1; then the weights $w_1, w_2, w_3$
of the one-dimensional representations satisfy
$w_1 + w_2 + w_3 = r$. The age of the action of a non-trivial group
element $e_i$ is given by
$\age(e_i, \C^3) =
\age(e_i^{w_1}) + \age(e_i^{w_2}) + \age(e_i^{w_3})$.  In this section
we will develop recursions for invariants of the form
\begin{equation} 				\label{eq:C3mr-invs}
\langle h_{e_1} \otimes \dots \otimes h_{e_{n-1}} \otimes \psi_n^\nu h_{e_n}
\rangle_{0, n}^{[\C^3/\mu_r]}
\end{equation}
where $e_1, \dots e_{n-1} \in \mu_r$ are group elements of age 1, and
$e_n$ is arbitrary, and also the only element for which we allow insertion
of a $\psi$-class. (This implies that
$\delta^{(1)}_T + \delta^{(2)}_T + \delta^{(3)}_T =
\abs{T}$ for all $T \subset [n-1]$.)

We want to determine the integral of a summand on the right-hand-side of
equation (\ref{eq:incl-excl-applied}) related to a partition
$\PP$ of $[n-1]$, after inserting an additional
$\psi$-class at the $n$-th marked point. Let $T_1, \dots, T_k$ be
the elements of $\PP_{\ge 2}$.
\[
\psi_n^\nu \cdot c(\PP) = \psi_n^\nu \cdot
\prod_{i=1}^k x_{T_i}
\prod_{\substack{T \subsetneq [n-1]/(\PP) \\ 2 \le \abs{T}}} (1 + x_T)
\prod_{a=1}^N
\prod_{p= \fract{\delta^{(a)}_{[n-1]}}}^{\delta^{(a)}_{[n-1]} - 1}
(t_a - p \psi_n)
\]
This term is supported on the intersection
$D^{T_1} \cap \dots \cap D^{T_k}$, isomorphic to
$\Mbar_{0, T_1 \cup \{*\}} \times \dots \times \Mbar_{0, T_k \cup \{*\}}
\times \Mbar_{0, [n]/(\PP)}$; so in order to determine the integral of
$c(\PP)\cdot \psi_n^\nu$, we will write it as a product of
$D^{T_1} \cdots D^{T_k}$ with factors that are pulled back
from one of the components above.

The numerator of $x_{T_j}$ is the only factor that has terms
coming from
$\Mbar_{0, T_j \cup \{*\}}$, while its denominator involves $\psi_{T_j}$, which is the $\psi$-class of the node corresponding
to the marking $T_j$ on $\Mbar_{0, [n]/(\PP)}$.
To examine the numerator more closely, we first factor out $D^{T_j}$:
\begin{align*}
\prod_{a,p} (1-t_a^{-1}p\psi_{T_j^C}) - \prod_{a,p} (1+t_a^{-1}p\psi_{T_j}) &= \sum_{k>0}
\beta_k ((-1)^k\psi_{T_j^C}^k - \psi_{T_j}^k) \\
&= D^{T_j}\cdot\sum_{k>0} \beta_k \sum_{\ell=0}^{k-1} (-\psi_{T_j^C})^{\ell}\psi_{T_j}^{k-1-\ell}.
\end{align*}

The largest power of $\psi_{T_j^C}$ which appears in the last expression
is $$\sum_{a=1}^3 (\delta^{(a)}_{T_j} - \fract{\delta^{(a)}_{T_j}}) - 1
= \abs{T_j} - 1 - \sum_{a=1}^3 \fract{\delta^{(a)}_{T_j}}.$$ As the dimension
of $\Mbar_{0, T_j \cup \{*\}}$ is $\abs{T_j} - 2$, the expression only
has a term in the top degree if
$\sum_{a=1}^3 \fract{\delta^{(a)}_{T_j}} = 1$, which means that $e_{T_j} = \prod_{i \in T_j} e_i$ acts with age $1$
on $\C^3$ and acts nontrivially in each coordinate direction.
By the balancing condition, $e_{T_j}$ is prescribing
the monodromy of the node as seen from the component corresponding to
$[n]/(\PP)$.
As the integral of $\psi_{T_j^C}^{\abs{T_j}-2}$ is one, the integral of
the above product on $\Mbar_{0, T_j \cup \{* \}}$ is
$\prod_{a=1}^3 t_a^{-\ceil{\delta^{(a)}_{T_j}}} (\delta^{(a)}_{T_j}-1)!$
if the condition on $e_{T_j}$ is satisfied, and 0 otherwise.

On $\Mbar_{0, [n]/(\PP)}$, we are left with the following product:
\[
\psi_n^\nu \cdot
\prod_{a=1}^3 t_a^{-1 + \delta^{(a)}_{[n]}}
\prod_{\emptyset \neq T \subseteq \otn{n-1}/(\PP)}
\prod_{p = \fract{\delta^{(a)}_T}}^{\delta^{(a)}_T - 1}
	\frac{t_a + p \psi_T + p D^T}{t_a + p \psi_T}
\]
Here we used the same conventions as for formula
(\ref{eq:chern-class-compact}), applied to the set $[n-1]/(\PP)$
(so for example $D^{\{T_j\}}$ is identified with $- \psi_{\{T_j\}}$,
which is the $\psi$-class \emph{of a single marking}; the term related
to $T_j$ in the above product is the denominator of $x_{T_j}$); and we
extended $\delta^{(a)}_T$ in the obvious way from subsets of $[n-1]$ to
subsets of the quotient set $[n-1]/(\PP)$.
Now by the remarks in section \ref{subsect:formula-remarks},
based on lemma \ref{combinatorial_relation}, this product
computes the Chern class of the obstruction bundle
on $\Mbar_{0, [n]/(\PP)} \( (e_i)_{i \in [n]/(\PP)}; B \mu_r\)$.
Its integral is thus given by the equivariant Gromov-Witten invariant
\[
\langle \bigotimes_{i \in [n-1]/(\PP) } h_{e_i} \otimes \psi_n^\nu h_{e_n}
	\rangle_{0, [n]/(\PP)}^{[\C^3/\mu_r]}.
\]

This proves the following  recursion:
\begin{Prop}		\label{prop:non-equiv-recursion}
For an equivariant Gromov-Witten invariant of $[\C^3/\mu_r]$ as in
equation (\ref{eq:C3mr-invs}) with the
assumptions above, let $S$ be the set of non-trivial partitions
$\PP$ of $[n-1]$
such that for every $T \in \PP_{\ge 2}$, the group element
$e_{T_j} = \prod_{i \in T_j} e_i$ acts with age 1 on $\C^3$, and non-trivial
in every coordinate direction.
Then the following recursive formula holds:
\begin{multline*}
\langle h_{e_1} \otimes \dots \otimes h_{e_{n-1}} \otimes
\psi_n^\nu h_{e_n} \rangle_{0, n}^{[\C^3/\mu_r]}
=
\langle h_{e_1} \otimes \dots \otimes h_{e_{n-1}} \otimes
\psi_n^\nu h_{e_n} \rangle_{0, n}^{[\C^3/\mu_r]; \text{weighted}} +
\\
+ \sum_{\PP \in S} (-1)^{\abs{\PP_{\ge 2}}+1}
  \prod_{T \in \PP_{\ge 2}} \prod_{a=1}^3 (\delta_T^{(a)}-1)!
\ \langle \bigotimes_{i \in [n-1]/(\PP) } h_{e_i} \otimes \psi_n^\nu h_{e_n}
	\rangle_{0, [n]/(\PP)}^{[\C^3/\mu_r]}.
\end{multline*}
\end{Prop}
Here the invariant with superscript ``weighted'' means the invariant as
computed by using the moduli space of weighted stable maps instead of the
ordinary moduli space, with weight data chosen as in section
\ref{subsect:Pn-3-computation}. These invariants are given, up to a
multiplication with a monomial in the $t_a$, by the
$(n-3-\nu)$-th elementary symmetric function of the variables
$t_a^{-1} (\delta^{(a)}_{[n-1]} - 1), t_a^{-1} (\delta^{(a)}_{[n-1]} - 2),
\dots, t_a^{-1} \fract{\delta^{(a)}_{[n-1]}}$ for $a = 1,2, 3$.

\subsection{Recursions for $[\C^3/\mu_3]$}
The recursion of proposition \ref{prop:non-equiv-recursion} simplifies further
in the case $[\C^3/\mu_3]$ for the diagonal representation of $\mu_3$.
The only group element
of age 1 is $\omega = e^{\frac{2 \pi i}3}$. We have
$\delta_T^{(a)} = \frac{\abs{T}}3$ for all $T \subset [n-1]$ and $a = 1,2,3$.
The set $S$ contains the
partitions $\PP$ of $[n-1]$ so that every $T \in \PP$ has
size $3 m_T + 1$ for some $m_T \in \Z_{\ge 0}$. The summand for $\PP$ in
the formula of proposition \ref{prop:non-equiv-recursion} depends only
on the sizes of the subsets, not on the actual subsets; if we set
$p = \sum_{T} m_T$, we can thus reduce the above sum to a sum over
partitions $\m = (m_1, \dots, m_k)$ of $p$, for all $p \ge 1$ with
$n - 3p \ge 3$. For any such partition, let
$M (n-1, \m)$ be the multinomial coefficient
\[
M(n-1, \m) =
	\binom{n-1}{3m_1 + 1, \dots, 3 m_k + 1, n-1-\sum_j (3m_j + 1)}
\]
counting the ways to distribute $n-1$ markings on the different
components.

\begin{Prop}			\label{prop:C3Z3-recursion}
\begin{multline*}
\langle h_{\omega}^{\otimes n-1} \otimes \psi_n^\nu h_{e_n}
\rangle_{0,n}^{[\C^3/\mu_3]}
=
\langle h_{\omega}^{\otimes n-1} \otimes
\psi_n^\nu h_{e_n} \rangle_{0, n}^{[\C^3/\mu_3]; \text{weighted}} \\
+ \sum_{p=1}^{\floor{\frac{n-3}3}} \sum_{\m} \frac{(-1)^{\abs{\m}+1}}{|\Aut \m|}
	\prod_{j=1}^k \((m_j - \frac 23)!\)^3
	M(n-1, \m)
	\langle h_{\omega}^{\otimes n-1 - 3p} \otimes \psi_n^\nu h_{e_n}
	\rangle_{0,n}^{[\C^3/\mu_3]}
\end{multline*}
\end{Prop}
The superscript ``weighted'' means the same as before; if $e_n = \omega$
and $\nu = 0$, this weighted Gromov-Witten invariant is just given as
$(-1)^{n+1}\((\frac{n-4}{3})!\)^3\frac{1}{3}$; otherwise it is an elementary symmetric function.

\begin{figure}[htb] \label{fig:recursion}
\includegraphics{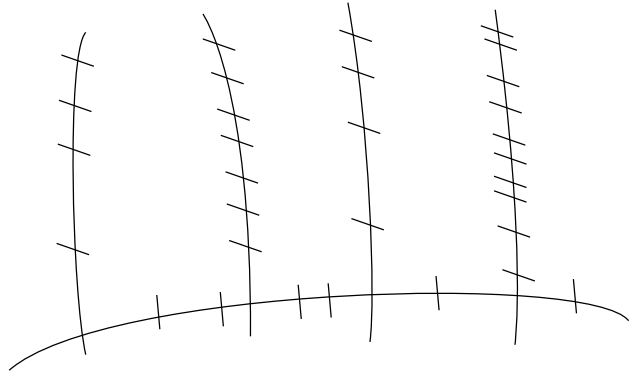}
\caption{Comb for $n = 30$, $p= 7$ and $\m = (1, 1, 2, 3)$}
\end{figure}

A maple program implementing some of these recursions is available from the
authors upon request. The numbers match the calculations of
\cite{ABK:topological_strings}, \cite{CCIT:computing} and \cite{Chuck-Renzo}.

As the recursions are linear, it is not hard to invert the matrix and obtain a direct formula.  Let $I_{\ell}=\langle h_{\omega}^{\otimes 3\ell+3}\rangle_{0,3\ell+3}^{[\C^3/\mu_3]}$.  Proposition \ref{prop:C3Z3-recursion} implies that
$$\sum_{p=0}^\ell (-1)^p C_{p,\ell}I_{\ell-p} = (-1)^\ell((\ell-\frac 13)!)^3\frac 13,$$
where $C_{0,\ell}=1$ and for $0<p\le \ell$, $C_{p,\ell}$ is the sum over partitions $\m=(m_1,\ldots,m_k)$ of $p$, with $k\le 3(\ell-p)+2$, of the quantity
$$\frac{1}{|\Aut \m|}\prod_{j=1}^{k} ((m_j- \frac 23)!)^3\cdot \binom{3\ell+2}{3m_1+1,\ldots,3m_k+1, 3(\ell-p)+2-k}.$$
Let $D_{p,\ell} = C_{\ell-p,\ell}$.  By inverting the matrix we obtain the formula:
$$3(-1)^\ell I_\ell = \sum_{\substack{S\subseteq [0,\ell-1] \\ S = \{x_0,\ldots,x_q\}}} (-1)^{|S|} ((x_0-\frac 13)!)^3 D_{x_0,x_1}\cdots D_{x_{q-1},x_q}D_{x_q,\ell},$$
where $[0,\ell-1]=\{0,1,\ldots,\ell-1\}$ and it is assumed that $x_0<x_1<\cdots<x_q$.  For $S=\emptyset$, the summand is taken to be $((\ell-1/3)!)^3$.

\subsection{Inversion of the ``mirror map''}
\label{recursion-and-mirror-map}
The recursion in this case can also be derived from the results in
\cite{CCIT:computing}. We explained above in section \ref{twisted-I}
that their $I$-function $I^{\mathrm{tw}}(t)$
is almost identical (up to high powers of $z$)
to the ``weighted $J$-function'' $J^{X; \text{weighted}}$.

  From general principles of Givental's formalism they deduce that for
the coordinate change $\tau(t) \colon H \to H$, called the
``mirror map'', given by
\[
\tau(t_0 h_1 + t_1 h_{\omega}) = t_0 h_{1}
+ \sum_{k \ge 0} \frac{(-1)^{3k} (t_1)^{3k+1}}{(3k+1)!}
	\bigl((k-\frac 23)!\bigr)^3 h_{\omega}
\]
the twisted $I$-function and the twisted $J$-function can be related
after setting the dual coordinate of $h_{\omega^2}$ equal to zero:
\[
I^{\mathrm{tw}}(t_0 h_1 + t_1 h_{\omega}, z)
= J^X (\tau(t_0 h_1 + t_1 h_{\omega}), z)
\]
Our recursive formula can be recovered by comparing coefficients of these
two power series; in other words, the sum over partitions $\PP$ in
proposition \ref{prop:C3Z3-recursion} is a combinatorial inversion
of the mirror map $\tau(t)$.

More precisely, given any two power series $A(t_0 h_1 + t_1 h_\omega, z)$,
$B(t_0 h_1 + t_1 h_\omega, z)$ related by
\[A(t_0 h_1 + t_1 h_\omega, z)
= B(\tau(t_0 h_1 + t_1 h_\omega), z)\]
we can recursively recover the coefficients
of $B$ by comparing coefficients of powers of $t_1$. If we write
$A(t_1 h_1) = \sum_{k} \frac{a_k}{k!} t_1^k$ and
$B(t_1 h_1) = \sum_{k} \frac{b_k}{k!} t_1^k$ with
$a_k, b_k \in H[[z^{-1}]]$, the recursion
will look exactly as proposition \ref{prop:C3Z3-recursion} with
the ``weighted invariants`` replaced by $a_k$ and the actual invariants
replaced by $b_k$. In particular, setting $A$ to the identity power series
$A(t_0 h_1 + t_1 h_\omega) = t_0 h_1 + t_1 h_\omega$
yields an inversion of the mirror map that can
also be a interpreted as a sum over comb curves.

\subsection{Equivariant recursions}
The methods of this section are sufficient to produce a linear recursion for
the equivariant descendant Gromov-Witten invariants of $[\C^N/\mu_r]$.
However, this requires one to allow $\psi$-classes at every marked point.
As in section \ref{s:non_eq_c3mur}, one can use lemma
\ref{lem:inclusion-exclusion}
to expand the equivariant Euler class of the obstruction bundle, and to each
partition of $[n-1]$ one should associate a comb as before,
where the $n$-th marked point is on the
head of the comb.  For each tooth of the comb, we can write the
numerator of $x_T$ as
\[
D^T\cdot \sum_{k>0} \beta_k \sum_{\ell=0}^{k-1} (-\psi_{T^C})^{\ell}\psi_T^{k-1-\ell},\] just as in section \ref{s:non_eq_c3mur}.
The exponent of $\psi_{T^C}$ which leads to a nonzero integral is determined
by the descendant exponents chosen for the marked points in $T$.  The integral
over the tooth can then be computed using the well-known formula
\[ \int_{\overline{M}_{0,n}}\psi_1^{a_1}\cdots\psi_n^{a_n} =
\binom{n-3}{a_1,\ldots,a_n}.\]
We are left with a polynomial in $\psi_T$, the $\psi$-class of the node
on the main component (and hence our method does not yield a recursion
for non-descendant invariants only).

In summary, for each partition of $[n-1]$, one gets a linear
combination of equivariant Gromov-Witten invariants with fewer marked points,
each with a combinatorial factor.  These must be summed together and added to
the weighted invariant, just as in section \ref{s:non_eq_c3mur}.

\appendix

\section{Combinatorics of divisors on $\Mbar_{0,n}$}

\subsection{Notations for divisors}

This section reviews notations for divisors on $\Mbar_{0,n}$, some of which
is introduced in this paper.  The standard relations are reviewed, and a
combinatorial proof of a key simplifying relation, used in section
\ref{sect:recursions}, is worked out.

First recall the vital divisor $D^T$, introduced by Keel \cite{Keel:M0n}, where $T$ is any subset of $[n]$ having at least $2$ and at most $n-2$ elements.  This
divisor is the locus of curves having a node which separates the markings into
$T$ and $T^C$.  Here complements are always taken within $[n]$.  To make
this into a correspondence, assume throughout this appendix that $S$ and $T$ are subsets of $[n-1]$.  It is natural to define $D^{[n-1]}$ to be $-\psi_n$, which comes from the work of de Concini and Procesi \cite{deConPro:hyperplane}.  From their point of view,
$D^{[n-1]}$ is the pullback of minus the hyperplane class under a sequence of
blowups producing $\Mbar_{0,n}$ from $\P^{n-3}$.  This sequence of blowups is the same one discussed in section \ref{sect:Pn-weights} in the context of weighted stable maps.
Under this blowup description of $\Mbar_{0,n}$, $D^T$ is the exceptional divisor of the blowup in the proper transform of the linear
space generated by the points labeled by $[n-1]\setminus T$.

The ring $H^*(\Mbar_{0,n})$ is generated by the divisors
$D^T$ for $T\subseteq[n-1]$, $|T|\ge 2$, with relations given by
\begin{align}
D^SD^T &\;\; \text{if $S$ and $T$ are incomparable and $S\cap T\not=\emptyset$},
\label{deConProc1}\\
\sum_{i,j\in T} D^T &\;\; \text{for every $i\not=j\in [n-1]$}\; \text{\cite{deConPro:hyperplane}}. \label{deConProc2}
\end{align}
Geometrically, the first relation is due to the fact that the exceptional
divisors $D^S$ and $D^T$ are disjoint, and the second is due to the fact
that the preimage of the hyperplane in $\P^{n-3}$ generated by all points
except $i$ and $j$ consists of the proper transform $D^{i,j}$ together with
all the exceptional divisors $D^T$ for $i,j\in T$.  Fix $T\subsetneq [n-1]$
and choose $i\notin T$ and $j\in T$.  Then relations
(\ref{deConProc1},\ref{deConProc2}) imply that
\begin{equation}
\label{deConProc3}
D^T\sum_{\substack{S:\; i\in S \\ S\supset T}} D^S = 0.
\end{equation}

For any nonempty subset $T\subseteq [n-1]$, we introduce the notation
\begin{equation} 			\label{eq:psi_T}
\psi_T :=\sum_{S\supsetneq T} D^S.
\end{equation}
If $T=\{i\}$, then this recovers the
$\psi$-class at the $i$-th marked point.  Indeed, we have for any distinct
$j,k$ different from $i$ (recalling the convention $n\notin T$),
$$\psi_i=\sum_{\substack{j,k\in T \\ i\notin T}} D^T + \sum_{\substack{i\in
T \\ j,k\notin T}} D^T = -\sum_{i,j,k\in T} D^T + \sum_{\substack{i\in T\\
j,k\notin T}} D^T + \sum_{i,j\in T} D^T + \sum_{i,k\in T} D^T =
\psi_{\{i\}},$$ the last equality following from an inclusion-exclusion
argument.  By definition, $\psi_{[n-1]} = 0$, and we also set
$D^{i}=-\psi_i$ for $1\le i\le n-1$.  Finally, we define $$\psi_{T^C}:=-D^T-\psi_T$$ for $T\subseteq[n-1]$, which we found useful in section \ref{sect:recursions}.

The Chern class formula of theorem \ref{mainthm} can now be expressed as
\begin{equation} \label{Chern_class_appendix}
c=\prod_{\emptyset\not= T\subseteq [n-1]}
\prod_{p=\fract{\delta_T}}^{\delta_T-1}
\left(\frac{1+p(D^T+\psi_T)}{1+p\psi_T}\right),
\end{equation}
where $\delta_T=\sum_{i\in T}\delta_i$ and $\delta_i$ are chosen so that
$0\le\delta_i <1$.  The following lemma shows that if every $\delta_i\ge 0$, this
expression is periodic in each $\delta_i$ with period $1$; so it
defines a continuous, piecewise-analytic function from an $(n-1)$-dimensional real
torus into the cohomology of $\Mbar_{0,n}$ with real coefficients.  This is used in
section \ref{sect:recursions} to produce recursions for the Gromov-Witten invariants.

\begin{Lem}
\label{combinatorial_relation}
Let $\delta_1,\ldots,\delta_{n-1}$ be real numbers, and for any subset
$T\subseteq[n-1]$, define $\delta_T=\sum_{i\in T}\delta_i$.  Then
\begin{equation} \label{combinatorial_relation_formula}
\prod_{1\in T}\frac{1+\delta_T(D^T+\psi_T)}{1+\delta_T\psi_T}=1.
\end{equation}
\end{Lem}

\begin{proof}
For $1\le k\le n-1$, let
$$E_k =
\sum_{i=2}^{n-1}\delta_i\sum_{\substack{1,i\in T \\ |T|>k}} D^T +
\delta_1\sum_{\substack{1\in T \\|T|>k}} D^T,$$
 and let
$$A_k =
(1+E_k)\prod_{\substack{1\in T\subseteq[n-1] \\ |T|\le k}}
\frac{1+\delta_T(D^T+\psi_T)}{1+\delta_T\psi_T}.$$
Then it must be shown
that $A_{n-1}=1$.  As $E_{n-2}=\delta_{[n-1]}D^{[n-1]}$, it follows that
$A_{n-2}=A_{n-1}$.  Moreover, $E_1=\delta_1\psi_1$, so $A_1=1$.  It remains
to show that $A_k=A_{k-1}$ for $1<k<n-1$.

Note that for any divisors $x,y,z$,
\begin{equation}
\label{divisor_relation}
\frac{(1+x)(1+y)}{1+z}=1+x+y-z \;\; \text{if}\;\; (x-z)(y-z)=0.
\end{equation}
Fix $T$ with $1\in T$ and $|T|=k$, and let $x$ be any expression of the form
$E_k+\sum_{S\in\sigma} \alpha_SD^S$, where $\sigma$ is a collection of subsets
$S\subseteq[n-1]$ with $|S|=k$ and $1\in S\not= T$.  Then $D^T(x-\delta_T\psi_T)=0$
by the following argument.  For any $S\in\sigma$, $1\in S\cap T$ and $|S|=|T|$, so
$D^SD^T=0$.  Moreover, $$E_k-\delta_T\psi_T = \sum_{i:\; i\notin T}
\delta_i\sum_{\substack{S:\; i\in S \\ T\subseteq S}} D^S +
\sum_{\substack{|S|>k \\ S\not\supset T \\1\in S}} \delta_SD^S.$$  Since $D^T$
annihilates the first term on the right hand side by (\ref{deConProc3}) and
annihilates the second by (\ref{deConProc1}), it follows that
$D^T(E_k-\delta_T\psi_T)=0$, completing the verification of the claim.

It follows that one can iterate through all sets $T$ with $|T|=k$ to
eliminate those factors from the expression for $A_k$ and apply relation
(\ref{divisor_relation}) at each step.  Since $$E_k+\sum_{\substack{1\in T
\\ |T|=k}} \delta_TD^T = E_{k-1},$$ it follows that $A_k=A_{k-1}$,
completing the proof.  \end{proof}

If in (\ref{Chern_class_appendix}), one were to add $1$ to $\delta_1$ (and
thus add one to each $\delta_T$ for $1\in T$) then one would be multiplying
$c$ by the left hand side of (\ref{combinatorial_relation_formula}).
Therefore $c$ does not change after the translation $\delta_i\mapsto
\delta_i+1$.  Using the notational conventions of section
\ref{sect:notations}, the formula for $c$ makes sense for negative values of
$\delta_i$, and the same argument shows that it remains invariant under
integer translations.  It seems natural to regard $\delta_1,\ldots,\delta_n$
as coordinates on $(\R/\Z)^n$ satisfying $\sum_{i=1}^n\delta_i=0$.

\subsection{Restricting to $D^T$}

It is a standard fact for any proper subset $T\subset [n-1]$ containing at
least two elements, $D^T\cong \Mbar_{0,|T|+1}\times\Mbar_{0,n-|T|+1}$, with
the node counting as an extra marked point on each factor.  The
restrictions of divisors $D^S$ and $\psi_S$ to $D^T$ are easily computed if
one uses subsets of $T$ for the divisors on the first factor on subsets of
the quotient set
$[n-1]/T$ for divisors on the second factor.
So on the first factor, the node counts as the extra
marked point that all sets must avoid.  For sets $S\subset T$, we use $D^S_*$ and $\psi^*_S$ to notate divisors on $\Mbar_{0,|T|+1}$ and we likewise use
$D^S_{\bullet}$ and $\psi_S^{\bullet}$ for divisors on $\Mbar_{0,n-|T|+1}$.
Now the following formulas hold for
any nonempty set $S\subseteq [n-1]$.
\begin{gather}
D^S|_{D^T} =
\begin{cases}
D^S_*\otimes 1, &\text{if $S\subsetneq T$} \\
1\otimes D^S_{\bullet}, &\text{if $S\cap T=\emptyset$} \\
1\otimes D^{S/T}_{\bullet}, &\text{if $S\supsetneq T$} \\
D^T_*\otimes 1 + 1\otimes D^{T/T}_{\bullet}, &\text{if $S=T$} \\
0, &\text{otherwise}
\end{cases} \label{restriction_Ds_Dt} \\
\psi_S|_{D^T} =
\begin{cases}
\psi_S^*\otimes 1, &\text{if $S\subsetneq T$} \\
1\otimes \psi_S^{\bullet}, &\text{if $S\cap T=\emptyset$} \\
1\otimes \psi_{S/T}^{\bullet}, &\text{if $S\supseteq T$}
\end{cases} \label{restriction_Ps_Dt}
\end{gather}
Note that the fourth line of (\ref{restriction_Ds_Dt}) is another way of writing the standard fact that the restriction of $-D^T$ to $D^T$ is the sum of the $\psi$-classes at the node on the two components.  Moreover, the third line of (\ref{restriction_Ps_Dt}) shows that $\psi_T$ restricts to $D^T$ as the $\psi$ class of the node on the component corresponding to $T^C$.  So our definition of $\psi_T^C$ as $-D^T-\psi_T$ ensures that $\psi_T^C$ restricts to the $\psi$ class of the node on the other component.

\addcontentsline{toc}{section}{References}
\bibliography{all}                      
\bibliographystyle{halpha}     

\end{document}